%% file: tempered.tex
\documentclass[11pt]{report}
\usepackage{amsmath}
\usepackage{amssymb}
\usepackage{amsthm}
\usepackage{latexsym}
\usepackage{graphicx}
\usepackage{hyperref}

\newtheorem{thm}{Theorem}[chapter]
\newtheorem{cor}[thm]{Corollary}
\newtheorem{lem}[thm]{Lemma}
\newtheorem{prop}[thm]{Proposition}

\providecommand{\norm}[1]{\left\| #1 \right\|}

\newcommand{\mh}{\mathbb}
\newcommand{\mr}{\mathrm}
\newcommand{\mc}{\mathcal}

\newcommand{\ds}{\displaystyle}
\newcommand{\scs}{\scriptstyle}

\newcommand{\ep}{\epsilon}
\newcommand{\es}{\emptyset}
\newcommand{\af}{\mr{aff}}
\newcommand{\hot}{\widehat{\otimes}}

\newcommand{\wig}{\textstyle \bigwedge}
\newcommand{\inp}[2]{\langle #1 \,,\, #2 \rangle}

\begin{document}

\include{tempered0}

\include{tempered1}

\include{tempered2}

\include{tempered3}

\include{temperedA}

\end{document}

%% file: tempered0.tex
\ \vspace{1cm}
\begin{center}
\huge \textbf{Homological algebra for affine Hecke algebras}\\[1cm]
\Large Eric Opdam and Maarten Solleveld\\[3mm]
\normalsize Korteweg-de Vries Institute for Mathematics,
Universiteit van Amsterdam\\
Plantage Muidergracht 24, 1018TV Amsterdam, The Netherlands\\
email addresses: opdam and mslveld at science.uva.nl \\[5mm]
August 2007 \\[2cm]
\end{center}

\begin{minipage}{13cm}
\textbf{Mathematics Subject Classification (2000).} \\
20C08; 18Gxx; 20F55 \\[1cm]
\textbf{Abstract.}\\
In this paper we study homological properties of modules over an
affine Hecke algebra $\mc{H}$. In particular we prove a comparison
result for higher extensions of tempered modules when passing to
the Schwartz algebra $\mc{S}$, a certain topological completion of
the affine Hecke algebra. The proof is self-contained and based on
a direct construction of a bounded contraction of certain standard
resolutions of $\mc{H}$-modules.

This construction applies for all positive parameters of the
affine Hecke algebra. This is an important feature since it is an
ingredient to analyse how the irreducible discrete series
representations of $\mc{H}$ arise in generic families over the
parameter space of $\mc{H}$. For irreducible non-simply laced
affine Hecke algebras this will enable us to give a complete
classification of the discrete series characters for all positive
parameters (we will report on this application in a separate
article).
\\[1cm]
\textbf{Acknowledgements.}\\
This research originated from joint work of Mark Reeder and
the first author. We are also grateful to Ralf Meyer and Henk Pijls
for some useful advice.

A part of this paper was written while the authors were guests of the
Max Planck Institut f\"ur Mathematik and the Hausdorff Research
Institute for Mathematics, both in Bonn. We thank these institutions
for their hospitality.
\end{minipage}

\tableofcontents

\chapter*{Introduction}
\addcontentsline{toc}{chapter}{Introduction}

Affine Hecke algebras are useful tools in the study of the
representation theory and harmonic analysis of a reductive $p$-adic
group $G$, cf. \cite{BuKu1,BuKu2,Lus3,Mor1,Mor2}. A central theme in
this context is the Morita equivalence of Bernstein blocks of the
category of smooth representations of $G$ with the module category
of suitable Hecke algebras, often closely related to affine Hecke
algebras. This could be thought of as an
affine analogue of the role played by finite dimensional
Iwahori-Hecke algebras in the representation theory of finite
groups of Lie type, a theory which was developed in great detail
by Howlett and Lehrer \cite{HoLe}.
An important point of Howlett-Lehrer theory
is the fact that the Hecke algebras which arise are semisimple
specializations of a generic algebra. The affine Hecke algebras
which arise in the study of reductive $p$-adic groups are
specializations of generic algebras as well. This time however, it
is much more delicate to relate the representation theory of
different specializations of the generic algebra. The theory
developed in this paper gives an important handle on such problems.

Various aspects of the harmonic analysis on $G$ can be transferred
to Hecke algebras \cite{HeOp}. In particular the Hecke algebra comes
equipped with a Hilbert algebra structure defined by an anti-linear
involution and a tracial state whose spectral measure (also called
Plancherel measure) corresponds to the restriction of the Plancherel
measure of $G$ to the Bernstein block under the Morita equivalence.
This should be compared to the role of generic degrees of
representations of finite dimensional Hecke algebras in
Howlett-Lehrer theory.
\\[2mm]

The Schwartz algebra completion $\mc{S}$ of $\mc{H}$ plays a role
which is similar to that of the Harish-Chandra Schwartz space $\mc C
(G)$ in the representation theory of $G$. In particular the support
of the Plancherel measure of $\mc H$ consists precisely of the
irreducible representations which extend continuously to $\mc S$
(the irreducible tempered representations).

More restrictively we say that an irreducible $\mc H$-module
belongs to the discrete series if it is contained in the left
regular representation of $\mc H$ on its own Hilbert space
completion. Every irreducible representation can be constructed
from a discrete series representation, with a suitable version
of parabolic induction. Therefore the discrete series is of
utmost importance in the representation theory of $\mc H$ and of
$\mc S$.

Although $\mc S$ is larger then $\mc H$, its
representation theory is actually simpler. The spectrum of $\mc S$
(also called the tempered spectrum of $\mc H$) is much smaller
than the spectrum of $\mc H$. For example the discrete series
corresponds to isolated points in the spectrum of $\mc S$, while
the spectrum of $\mc H$ is connected. This observation leads to an especially
nice property of $\mc S$, namely that discrete series representations
are projective and injective as $\mc S$-modules. Contrarily $\mc H$
does not have finite dimensional projective modules. Yet with quite
some representation theory \cite{DeOp} one can reconstruct the
entire spectrum of $\mc H$ from its tempered spectrum.

A priori there could exist higher extensions of tempered
$\mc H$-modules which are themselves not tempered. But this does
never happen. More precisely we prove in Corollary \ref{cor:3.2} that
\begin{equation}\label{eq:1.1}
\mr{Ext}^n_{\mc H} (U,V) \cong \mr{Ext}^n_{\mc S} (U,V)
\end{equation}
for all finite dimensional tempered $\mc H$-modules $U$ and $V$ and
all $n \geq 0$. Our belief that something like \eqref{eq:1.1} might
be true was inspired by the work of Vign\'eras, Schneider, Stuhler
and Meyer \cite{Vig,ScSt,Mey2}.

To prove \eqref{eq:1.1} we construct explicit resolutions of $U$ and
$V$ by projective $\mc H$-modules. The remarkable part of the
proof is that we can turn these into projective $\mc S$-module
resolutions in the most naive way, simply by tensoring them with
$\mc S$ over $\mc H$.

One instance of \eqref{eq:1.1} is particularly important. Suppose
that $U$ is a discrete series representation and that $V$ is an
irreducible tempered $\mc H$-module. Theorem \ref{thm:3.5} states
that
\begin{equation}\label{eq:1.2}
\mr{Ext}_{\mc H}^n (U,V) \cong \left\{ \begin{array}{lll}
\mh C & \mr{if} \; U \cong V \; \mr{and} \; n = 0 \\
0 & \mr{otherwise}
\end{array} \right.
\end{equation}
We want to use \eqref{eq:1.2} to count the number of inequivalent
discrete series representations. This requires quite a few step,
which we discuss now. The Euler-Poincar\'e characteristic \cite{ScSt}
of two finite dimensional $\mc H$-modules is defined as
\begin{equation}\label{eq:1.3}
EP_{\mc H} (U,V) = \sum_{n=0}^\infty (-1)^n \dim_{\mh C}
\mr{Ext}^n_{\mc H} (U,V)
\end{equation}
This extends to a symmetric, bilinear and positive semidefinite
pairing on virtual $\mc H$-modules. By \eqref{eq:1.2} the discrete
series form an orthonormal set for this pairing.

On the other hand for the label function $q \equiv 1$ we have $\mc H
(\mc R ,1) = \mh C [W]$ and $\mc S (\mc R ,1) = \mc S (W)$, so
\eqref{eq:1.3} becomes
\begin{equation}\label{eq:1.4}
EP_W (U,V) = \sum_{n=0}^\infty (-1)^n \dim_{\mh C} \mr{Ext}^n_W
(U,V)
\end{equation}
This is much simpler than \eqref{eq:1.3}, as everything about the
Euler-Poincar\'e characteristic for groups like $W$ can be made
explicit. In Theorem \ref{thm:3.10} we find a conjugation-invariant
"elliptic" measure $\mu_{ell}$ on $W$ such that
\begin{equation}\label{eq:1.5}
EP_W (U,V) = \int_W \overline{\chi_U} \chi_V d \mu_{ell}
\end{equation}
where $\chi$ denotes the character of a representation. The support
of $\mu_{ell}$ consists precisely of the elliptic conjugacy
classes in $W$, whose number can easily be counted. This can
be compared with Kazhdan's elliptic integrals \cite{Kaz,ScSt,Bez}.

Finally we relate $EP_{\mc H}$ to $EP_W$ as follows. The label
function $q$ can be scaled to $q^\ep \; (\ep \in \mh R )$, which
yields a continuous field of algebras $\mc H (\mc R ,q^\ep )$. One
can associate to any finite dimensional $\mc H$-module $V$ a
continuous family of modules $\tilde \sigma_\ep (V)$ such that
\begin{equation}\label{eq:1.6}
EP_{\mc H (\mc R ,q^\ep )} \big( \tilde \sigma_\ep (U), \tilde
\sigma_\ep (V) \big) = EP_{\mc H}(U,V) \qquad \forall \ep \in [-1,1]
\end{equation}
In particular we can evaluate this at $\ep = 0$, which in
combination with the above yields a important upper bound on the
number of discrete series representations of $\mc H$, see
Proposition \ref{prop:3.9}. In \cite{OpSo} we will use this bound to
obtain a complete classification of the discrete series of affine
Hecke algebras $\mc H (\mc R ,q)$ with $\mc R$ irreducible and $q$
positive.
\\[2mm]

Now let us describe the contents of the chapters.

In the first
chapter we collect some notations and results that will be used
subsequently. We do not prove any deep theorems in this chapter, but
some of the results have not been published in research papers
before.

Chapter two is the technical heart of the paper, here we prove
everything needed for \eqref{eq:1.1}. In fact we do something
better, we construct an explicit projective $\mc H$-bimodule
resolution of $\mc H$. The crucial point is that this becomes a
resolution of $\mc S$ if we tensor it with $\mc S \otimes \mc
S^{op}$ over $\mc H \otimes \mc H^{op}$ and subsequently complete it
to a complex of Fr\'echet spaces. As an immediate consequence we
calculate that the global dimensions of $\mc H$ and $\mc S$ are
equal to the rank of the underlying root datum $\mc R$.

Although the proof of \eqref{eq:1.1} uses the combinatorial
structure of affine Hecke algebras in an essential way, the result
itself is of a more analytical nature. The inclusion $\mc H \to \mc
S$ can be compared to embeddings of the type $F_1 (G) \to F_2 (G)$,
where $G$ is a locally compact group and the $F_i (G)$ are certain
convolution algebras of functions on $G$. In many situations of this
type there is a comparison result
\begin{equation}
\mr{Ext}^*_{F_1 (G)} (U,V) = \mr{Ext}^*_{F_2 (G)} (U,V)
\end{equation}
for very general modules $U$ and $V$ \cite{Mey2}.

We choose to formulate our results in the category of bornological
$\mc S$-modules. Bornologies are the best technique to cover both
non-topological algebras like $\mc H$ and Fr\'echet algebras like
$\mc S$, in a natural way. However, we would like to point out that
the technical language of bornologies is inessential when dealing
with the case of finite dimensional modules of $\mc H$ or $\mc S$.
In this case it suffices to work with algebraic tensor products and
all proofs can be adapted in such a way so as to avoid the use of
results on bornologies. In particular the results on the discrete
series do not rely on bornologies. We have put some necessary
information on bornological modules in the appendix.

In chapter three we first study the Euler-Poincar\'e characteristic
for crossed products of lattices with finite groups. This leads
among others to \eqref{eq:1.5}. Clearly the results
hold for affine Weyl groups, but they do not rely on root systems.
In the last two sections we combine everything to derive the
aforementioned properties of the Euler-Poincar\'e characteristic
for affine Hecke algebras.

%% file: tempered1.tex
\chapter{Preliminaries}

\section{Root data}

First we introduce some well-known objects associated to root
data. For more background the reader is referred to
\cite{BrTi,Hum,IwMa}.

Let $R_0$ be a reduced root system of rank $r$ in an Euclidean
space $E \cong \mh R^r$. Let $W_0$ be the Weyl group of $R_0$ and
\[
F_0 = \{ \alpha_1 ,\ldots, \alpha_r \}
\]
an ordered basis. This determines the set of positive (resp.
negative) roots $R_0^+$ (resp. $R_0^-$). We suppose that $R_0$ is
part of a based root datum
\[
\mc R = (X, R_0, Y ,R_0^\vee ,F_0)
\]
For $I \subset F_0$ we write
\[
\begin{array}{lll}
C_I^+ & := & \{ x \in E : \inp{x}{\alpha_i^\vee} = 0 \;
\forall \alpha_i \in I ,\, \inp{x}{\alpha_j^\vee} \geq 0 \;
\forall \alpha_j \in F_0 \setminus I \}  \\
C_I^{++} & := & \{ x \in E : \inp{x}{\alpha_i^\vee} = 0 \;
\forall \alpha_i \in I ,\, \inp{x}{\alpha_j^\vee} > 0 \;
\forall \alpha_j \in F_0 \setminus I \}
\end{array}
\]
We call $C_\es^{++}$ the positive chamber. Its closure $C_\es^+$
is a fundamental domain for the action of $W_0$ on $E$.
The isotropy group (in $W_0$) of any point of $C_I^{++}$ is the
standard parabolic subgroup $W_I$ of $W_0$.

Recall that $Y \times \mh Z$ is the set of integral affine
linear functions on $X$. Let $R^\af$ be the affine root system
$R_0^\vee \times \mh Z \subset Y \times \mh Z$. The subsets of
positive and negative affine roots are
\[
\begin{array}{ccccc}
R_+^\af & = & R_0^{\vee,+} \times \{ 0 \} & \cup &
R_0^\vee \times \mh Z_{> 0} \\
R_-^\af & = & R_0^{\vee,-} \times \{ 0 \} & \cup &
R_0^\vee \times \mh Z_{< 0}
\end{array}
\]
The affine Weyl group of $R^\af$ is $W^\af = \mh Z R_0 \rtimes W_0$,
usually considered as a group of affine linear transformations of $X$.
It acts on $R^\af$ by
\[
w \cdot (\alpha^\vee ,k) (x) = (\alpha^\vee ,k) (w^{-1} x)
\]
For $a = (\alpha^\vee ,k) \in R^\af$ consider the affine hyperplane
\[
H_a := \{ x \in E : \inp{x}{a} = \inp{x}{\alpha^\vee} + k = 0 \}
\]
By definition $s_a$ is the reflection in this hyperplane, given
by the formula
\[
s_a (x) = x - \inp{x}{\alpha^\vee} \alpha - k \alpha
\]
Let $F_M$ be the set of maximal elements of $R_0^\vee$ for the dominance
ordering. Label its elements $\alpha_j^\vee ,\, j = r+1, \ldots, r+r'$,
where $r'$ is the number of irreducible components of $R_0$. We write
\[
a_j := \left\{ \begin{array}{ccc}
(\alpha_j^\vee ,0) & \mr{if} & \alpha_j^\vee \in F_0^\vee \\
(- \alpha_j^\vee ,1) & \mr{if} & \alpha_j^\vee \in F_M
\end{array} \right.
\]
Then
\[
F^\af := \{ a_j : j=1, \ldots ,r' \}
\]
is a basis of $R^\af$ and  $\left( W^\af ,S^\af \right)$
is a Coxeter system, where
\[
S^\af := \{ s_a : a \in F^\af \}
\]
For $J \subset S^\af$ we put
\[
A_J := \{ x \in E : \inp{x}{a_j} = 0 \: \forall a_j \in J ,\,
\inp{x}{a_i} > 0 \: \forall a_i \in F^\af \setminus J \}
\]
All the $A_J$ are facets of the fundamental alcove $A_\es$.
Its closure $\overline{A_\es}$ is a fundamental domain for
the action of $W^\af$ on $E$. \label{p:Sigma}
The isotropy group (in $W^\af$) of a point of $A_J$ is the
standard parabolic subgroup $\langle J \rangle$ of $W^\af$.
We will also write facets as $f = A_J$, in which case the
pointwise stabilizer is $W_f = \langle J \rangle$. Notice
that this is consistent with the above notation in the sense
that $W_0$ is the isotropy group of the facet $\{0\}$.

All the hyperplanes $H_{(\alpha^\vee ,k)}$ together give $E$
the structure of a polysimplicial complex $\Sigma$. The interior
of a polysimplex of maximal dimension is called an alcove.
\\[3mm]
\textbf{Example.}\\
Let $R_0$ be the root system $B_2$ in $E = \mh R^2$:
\[
R_0 = \{ \pm (1,-1), \pm (0,1), \pm (1,0), \pm (1,1) \}
\]
The Weyl group $W_0$ is isomorphic to the dihedral group
$D_4$. A basis of $R_0$ is
\[
F_0 = \{ \alpha_1 = (1,-1) , \alpha_2 = (0,1) \}
\]
The positive chamber and its walls are

\begin{picture}(6,6)(-1,-1)
\put(0,0){\line(1,1){4}}
\put(0,4){\line(1,-1){4}}
\put(2,0){\line(0,1){4}}
\put(0,2){\line(1,0){4}}
\put(1,1){\circle*{0.15}}
\put(2,1){\circle*{0.15}}
\put(3,1){\circle*{0.15}}
\put(3,2){\circle*{0.15}}
\put(3,3){\circle*{0.15}}
\put(2,3){\circle*{0.15}}
\put(1,3){\circle*{0.15}}
\put(1,2){\circle*{0.15}}
\put(2,2){\circle*{0.2}}
\put(3.1,1.1){\makebox{$\alpha_1$}}
\put(2.1,3.1){\makebox{$\alpha_2$}}
\put(4.1,1.8){\makebox{$C^{++}_{\{ \alpha_2 \} }$ }}
\put(4,4.1){\makebox{$C^{++}_{\{ \alpha_1 \} }$ }}
\put(3.7,2.7){\makebox{$C^{++}_{\es }$ }}
\newcounter{as}
\setcounter{as}{-3}
\multiput(-0.2,-0.4)(1,0){5}{\stepcounter{as}
  \makebox{\small{\arabic{as}}}}
\setcounter{as}{-3}
\multiput(-0.5,-0.1)(0,1){5}{\stepcounter{as}
  \makebox{\small{\arabic{as}}}}
\end{picture}
\\
If furthermore $\alpha_3 = (1,0)$ then
\[
F^\af = \{ (\alpha_1^\vee ,0) , (\alpha_2^\vee ,0) ,
(-\alpha_3^\vee ,1) \} = \{ a_1 ,a_2 ,a_0 \}
\]
The affine Weyl group $W^\af$ is generated by the simple reflections
\[
\begin{array}{ccccc}
s_1 & : & (x_1 ,x_2 ) & \to & (x_2 ,x_1 ) \\
s_2 & : & (x_1 ,x_2 ) & \to & (x_1 , -x_2 ) \\
s_0 & : & (x_1 ,x_2 ) & \to & (1 - x_1 ,x_2 ) \\
\end{array}
\]
The simplicial complex $\Sigma$ and the fundamental alcove look like

\begin{picture}(12,6)(-1,-1)
\multiput(0,0)(0.5,0){9}{\line(0,1){4}}
\multiput(0,0)(0,0.5){9}{\line(1,0){4}}
\put(0,0){\line(1,1){4}}
\put(0,4){\line(1,-1){4}}
\put(0,1){\line(1,1){3}}
\put(0,3){\line(1,-1){3}}
\put(0,2){\line(1,1){2}}
\put(0,2){\line(1,-1){2}}
\put(0,3){\line(1,1){1}}
\put(0,1){\line(1,-1){1}}
\put(1,0){\line(1,1){3}}
\put(1,4){\line(1,-1){3}}
\put(2,0){\line(1,1){2}}
\put(2,4){\line(1,-1){2}}
\put(3,0){\line(1,1){1}}
\put(3,4){\line(1,-1){1}}
\put(1,1){\circle*{0.15}}
\put(2,1){\circle*{0.15}}
\put(3,1){\circle*{0.15}}
\put(3,2){\circle*{0.15}}
\put(3,3){\circle*{0.15}}
\put(2,3){\circle*{0.15}}
\put(1,3){\circle*{0.15}}
\put(1,2){\circle*{0.15}}
\put(2,2){\circle*{0.2}}
\put(2.17,2.11){\makebox{$\scriptscriptstyle
  A_{\scriptscriptstyle \es}$}}
\setcounter{as}{-3}
\multiput(-0.2,-0.4)(1,0){5}{\stepcounter{as}
  \makebox{\small{\arabic{as}}}}
\setcounter{as}{-3}
\multiput(-0.5,-0.1)(0,1){5}{\stepcounter{as}
  \makebox{\small{\arabic{as}}}}

\put(6,1){\line(1,0){2}}
\put(6,1){\line(1,1){2}}
\put(8,1){\line(0,1){2}}
\put(6,1){\circle*{0.2}}
\put(8,1){\circle*{0.15}}
\put(8,1){\circle*{0.15}}
\put(8,3){\circle*{0.15}}
\put(7.2,1.5){\makebox{$A_\es$}}
\put(6.7,0.6){\makebox{$A_{\{a_2 \}}$}}
\put(8.1,2){\makebox{$A_{\{a_0 \}}$}}
\put(6.2,2.3){\makebox{$A_{\{a_1 \}}$}}
\put(8.1,0.8){\makebox{$A_{\{ a_0 ,a_2 \}} =\,
  \scs{\{ (1/2,0) \}}$}}
\put(8.1,3.1){\makebox{$A_{\{ a_0 ,a_1 \}} =\,
  \scs{\{ (1/2,1/2) \}}$}}
\put(5,0.7){\makebox{$A_{\{ a_1 ,a_2 \}}$}}
\put(5,0.2){\makebox{$=\, \scs{\{ (0,0) \}}$}}
\end{picture}

In general, if $A$ and $A'$ are two alcoves, then a gallery of
length $n$ between $A$ and $A'$ is a sequence $(A_0 ,\ldots ,A_n)$
of alcoves such that
\begin{itemize}
\item $A_0 = A$
\item $A_n = A'$
\item $\overline{A_{i-1}} \cap \overline{A_i} \:,\, \forall i$
is contained in exactly one hyperplane $H_a$
\end{itemize}
The group $W^\af$ acts simply transitively on the set of
alcoves. For $w \in W^\af$ there is a natural bijection
between expressions of $w$ in terms of the generators $S^\af$,
and galleries from $A_\es$ to $w A_\es$. This bijection is
given by
\begin{equation}\label{eq:2.1}
w = s_1 \cdots s_n \qquad \longleftrightarrow \qquad
(s_1 \cdots s_m A_\es )_{m=0}^n
\end{equation}

\begin{lem}\label{lem:2.1}
For $w \in W^\af$ the following numbers are equal:
\begin{description}
\item[1)] the word length $\ell (w)$ in the Coxeter system
$\left( W^\af ,S^\af \right)$
\item[2)] $\# \left\{ a \in R_+^\af : w a \in R_-^\af \right\}$
\item[3)] the number of hyperplanes $H_a \: (a \in R^\af )$
separating $A_\es$ and $w A_\es$
\item[4)] the minimal length of a gallery between $A_\es$ and $w A_\es$
\end{description}
In particular \eqref{eq:2.1} restricts to a bijection
between reduced expressions and galleries of minimal length.
\end{lem}
\emph{Proof.}
See \cite[Section 1]{IwMa}, \cite[Section 2.1]{BrTi} or
\cite[Theorem 4.5]{Hum}.
$\qquad \Box$\\[3mm]

Varying on the Bruhat order, we define a partial order $\leq_A$ on the
affine Weyl group $W^\af$:
\[
u \leq_A w \qquad \Longleftrightarrow \qquad
\ell (u) + \ell (u^{-1} w) = \ell (w)
\]
This means that $u \leq_A w$ if and only if a reduced expression for
$u$ can be extended to a reduced expression for $w$ by writing extra
terms on the right.

Let $K$ be a subset of $E$, and $\alpha \in R_0$.
\[
\begin{array}{ccr}
m(K,\alpha ) & := & \inf \left\{ \lfloor \inp{x}{\alpha^\vee} \rfloor :
x \in K \cup A_\es \right\} \\
M(K,\alpha ) & := & \sup \left\{ \lceil \inp{x}{\alpha^\vee} \rceil :
x \in K \cup A_\es \right\}
\end{array}
\]
where $\lfloor y \rfloor$ and $\lceil y \rceil$ denote respectively the
floor and the ceiling of a real number $y$. With these numbers we define
\[
\begin{array}{lll}
A(K,\alpha ) & := & \{ x \in E : m(K,\alpha ) \leq \inp{x}{\alpha^\vee}
\leq M(K,\alpha ) \} \\
A(K) & := & \bigcap_{\alpha \in R_0} A(K,\alpha )
\end{array}
\]
We can interpret $A(K)$ as a kind of $\Sigma$-approximation of the
convex closure of $K \cup A_\es$ in $E$.
\\[3mm]
\textbf{Example.}\\
In the setting of our previous example $R_0 = B_2$,
let $K$ be the simplex\\
$[(3/2,3/2),(3/2,2),(2,2)]$. Then $A(K)$ is the colored area below:
\\[1mm]

\includegraphics[width=77mm,height=4cm]{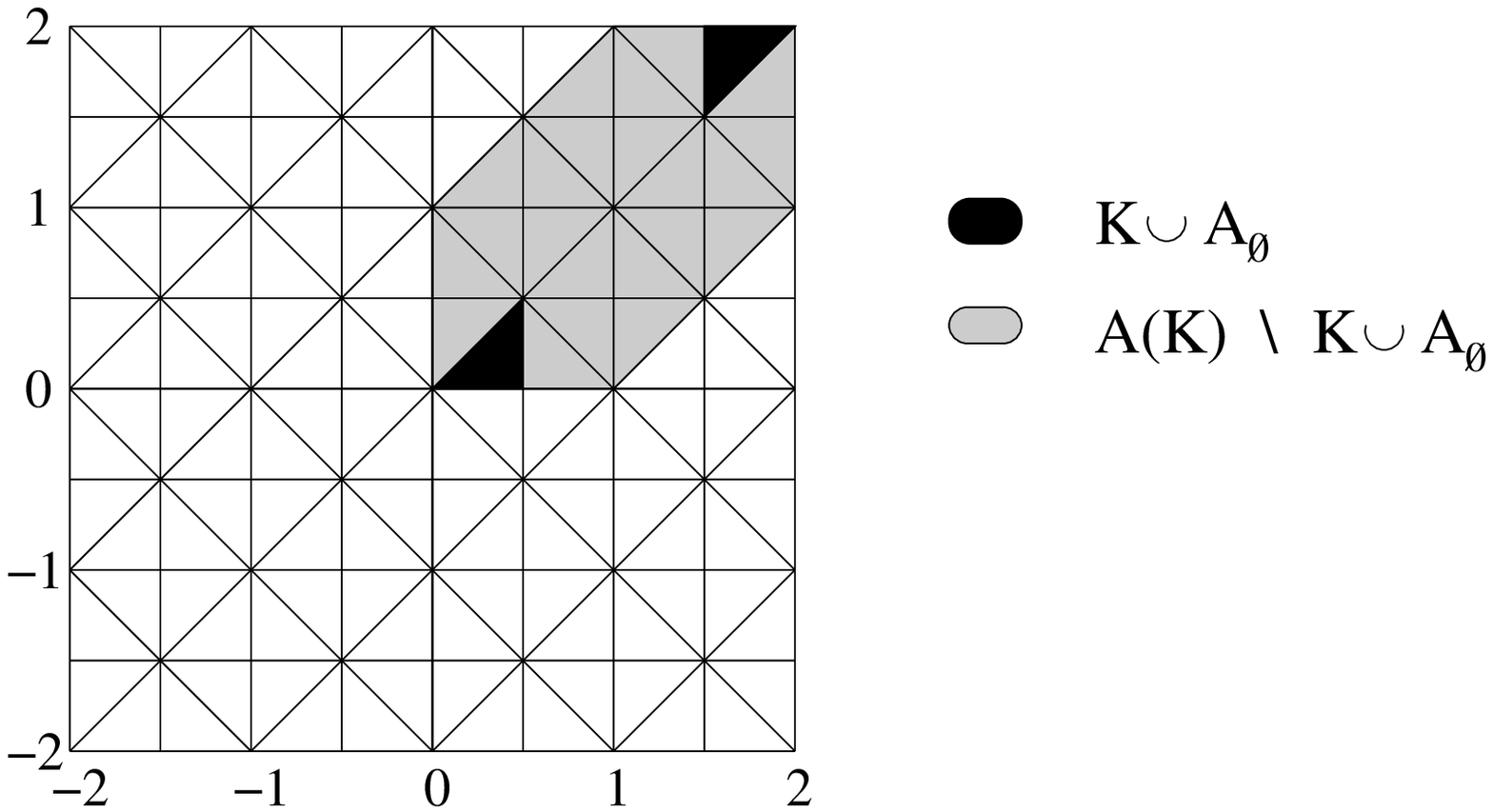}

\begin{lem}\label{lem:2.2}
For any $w \in W^\af$ we have
\[
A(w A_\es ) = \bigcup_{u \leq_A w} u \overline{A_\es}
\]
\end{lem}
\emph{Proof.} "$\supset$"
By Lemma \ref{lem:2.1} every alcove $u A_\es$ with $u \leq_A w$ is
part of a gallery of minimal length between $A_\es$ and $w A_\es$.
Such a gallery cannot cross any hyperplane
$H_a \; (a \in R^\af)$ that does not separate $A_\es$ and $w A_\es$.
So for every $\alpha \in R_0$ we must have
\[
\inp{u A_\es}{\alpha^\vee} \subset
[m(w A_\es ,\alpha ) ,M(w A_\es ,\alpha )]
\]
"$\subset$"
Since it is bounded by hyperplanes $H_a$ with
$a \in R^\af ,\, A(w A_\es )$ is a union of closures of alcoves.
If $B \subset A (w A_\es )$ is an alcove, then there are no
hyperplanes $H_a$ separating $B$ from $A_\es \cup w A_\es$.
Hence $B$ is part of at least one gallery of minimal length between
$A_\es$ and $w A_\es$. So $B = u A_\es$ for some $u \leq_A w.
\qquad \Box$
\\[3mm]
We note the consequence
\begin{equation}\label{eq:1.was}
w A (\sigma ) \subset A ( w \sigma ) \qquad
\forall \sigma \subset C_\es^+ , w \in W_0
\end{equation}
\newpage

\section{Affine Hecke algebras}

We recall a few important results on affine Hecke algebras,
meanwhile fixing some notations. Reconsider the based root datum
$\mc R = (X, R_0 ,Y, R_0^\vee ,F_0 )$. The Weyl group of $\mc R$ is
\[
W(\mc R ) = W = X \rtimes W_0
\]
which acts naturally on $X$. Clearly it contains $W^\af$ as
a normal subgroup. We write
\begin{align*}
&X^+ := \{ x \in X : \inp{x}{\alpha^\vee} \geq 0
\; \forall \alpha \in F_0 \}\\
&X^- := \{ x \in X : \inp{x}{\alpha^\vee} \leq 0
\; \forall \alpha \in F_0 \} = -X^+
\end{align*}
It is easily seen that the center of $W$ is the lattice
\[
Z(W) = X^+ \cap X^-
\]
We also want to make $W$ act on $E$. Since
\[
X \otimes \mh R = E \oplus \big( Z(W) \otimes \mh R \big)
\]
there is an orthogonal projection
\[
p_E : X \otimes \mh R \to E
\]
This induces a group homomorphism
\[
p_E : W \to E \rtimes W_0
\]
and the latter group acts naturally on $E$. The resulting action
of $W$ on $E$ consists of automorphisms of $\Sigma$, because
\[
\inp{p_E (x)}{\alpha^\vee} = \inp{x}{\alpha^\vee} \in \mh Z
\qquad \forall x \in X , \alpha^\vee \in R_0^\vee
\]
Hence 2), 3) and 4) of Lemma \ref{lem:2.1} define a natural
extension of the length function $\ell$ from $W^\af$ to $W$.

We say that $\mc R$ is semisimple if $R_0^\perp = 0 \subset Y$,
or equivalently if $X \otimes \mh R = E$. If $\mc R$ is not
semisimple then we can make it so by enlarging $R_0$ and
$R_0^\vee$. Namely, pick a basis
$\{ \alpha_{r+1} ,\ldots, \alpha_{\mr{rk}(X)} \}$ of
$X \cap (R_0^\vee )^\perp$ and declare those elements to be
simple roots. Also pick $\alpha_j^\vee \in Y$ such that
\[
\inp{\alpha_i}{\alpha_j^\vee} = 2 \delta_{ij} \qquad
i = 1, \ldots, \mr{rk}(X) ,\, j > r
\]
Thus we constructed a semisimple based root datum
\begin{equation}\label{eq:2.20}
\widetilde{\mc R} := (X, \widetilde R , 
Y, \widetilde R^\vee ,\widetilde F_0 )
\end{equation}
where $\widetilde R \cong R_0 \times (A_1 )^{\mr{rk}(X) - r}$.
Observe that
\begin{equation}\label{eq:2.29}
W (\widetilde{\mc R}) = W(\mc R ) \rtimes \tilde G =
X \rtimes \big( W_0 (\mc R ) \times \tilde G \big) =
X \rtimes W_0 (\widetilde{\mc R})
\end{equation}
With $\mc R$ we also associate some other root systems.
There is the non-reduced root system
\[
R_{nr} := R_0 \cup \{ 2 \alpha : \alpha^\vee \in 2 Y \}
\]
Obviously we put $(2 \alpha )^\vee = \alpha^\vee / 2$. Let $R_1$
be the reduced root system of long roots in $R_{nr}$:
\[
R_1 := \{ \alpha \in R_{nr} : \alpha^\vee \not\in 2 Y \}
\]
Let $q$ be a positive labeling of $R_{nr}^\vee$, i.e. a
$W_0$-invariant map $R_{nr}^\vee \to (0,\infty )$. This uniquely
determines a parameter function $q : W \to (0,\infty )$ with
the properties
\begin{equation}\label{eq:2.18}
\begin{array}{llll}
q (s_\alpha ) & = & q_{\alpha^\vee} & \alpha \in R_0 \cap R_1 \\
q (t_\beta s_\beta ) & = & q_{\beta^\vee} &
 \beta \in R_0 \setminus R_1 \\
q (s_\beta ) & = & q_{\beta^\vee /2} q_{\beta^\vee} &
\beta \in R_0 \setminus R_1 \\
q (\omega ) & = & 1 & \omega \in \Omega \\
q (w v) & = & q (w) q(v) & w,v \in W \quad
\mr{with} \quad \ell (wv) = \ell (w) + \ell (v)
\end{array}
\end{equation}
Conversely every function on $W$ with the last two properties
defines a labeling of $R_{nr}^\vee$. We speak of equal parameters if
$q(s) = q(s') \; \forall s,s' \in S^\af$.

The affine Hecke algebra $\mc H = \mc H (\mc R ,q)$ is the unique
complex associative algebra
with basis $\{ T_w : w \in W \}$ and relations
\[
\begin{array}{lllll}
T_w T_v & = & T_{wv} & \mr{if} & \ell (wv) = \ell (w) + \ell (v) \\
T_s T_s & = & (q(s) - 1) T_s + q(s) T_e & \mr{if} & s \in S^\af
\end{array}
\]
We can extend $q$ to a parameter function $\tilde q$ on
$W (\widetilde{\mc R} )$ by putting
\begin{equation}
\tilde q (s_{\alpha_j}) = 1 \quad \forall j > r
\end{equation}
By construction
\[
\mc H (\widetilde{\mc R},\tilde q ) \cong
\tilde G \ltimes \mc H (\mc R ,q)
\]
Now we describe the Bernstein presentation of $\mc H$.
For $x \in X^+$ we write
\[
\theta_x := N_x
\]
The corresponding semigroup morphism
$X^+ \to \mc H (\mc R ,q)^\times$
extends naturally to a group homomorphism
\[
X \to \mc H (\mc R ,q)^\times : x \to \theta_x
\]

\begin{thm}\label{thm:2.12}
\begin{description}
\item[a)] The sets $\{ T_w \theta_x : w \in W_0 , x \in X \}$ and
$\{ \theta_x T_w : w \in W_0 , x \in X \}$ are both bases of $\mc H$.
\item[b)] The subalgebra $\mc A := \mr{span} \{ \theta_x : x \in X \}$
is isomorphic to $\mh C [X]$.
\item[c)] The Weyl group $W_0$ acts on $\mc A$ by $w \cdot \theta_x =
\theta_{wx}$ and the center of $\mc H (\mc R ,q)$ is $Z (\mc H ) =
\mc A^{W_0}$.
\end{description}
\end{thm}
\emph{Proof.}
These results are due to Bernstein, see \cite[\S 3]{Lus2}. $\Box$
\newpage

Let $T$ be the complex algebraic torus $\mr{Hom}_{\mh Z} (X, \mh
C^\times )$, so that $\mc A \cong \mc O (T)$ and $Z (\mc H ) = \mc
A^{W_0} \cong \mc O (T / W_0 )$. From Theorem \ref{thm:2.12} we see
that $\mc H$ is of finite rank over its center, and hence
Noetherian.

For a set of simple roots $I \subset F_0$ we introduce the notations
\begin{equation}
\begin{array}{l@{\qquad}l}
R_I = \mh Q I \cap R_0 & R_I^\vee = \mh Q R_I^\vee \cap R_0^\vee \\
X_I = X \big/ \big( X \cap (I^\vee )^\perp \big) &
X^I = X / (X \cap \mh Q I ) \\
Y_I = Y \cap \mh Q I^\vee & Y^I = Y \cap I^\perp \\
T_I = \mr{Hom}_{\mh Z} (X_I, \mh C^\times ) &
T^I = \mr{Hom}_{\mh Z} (X^I, \mh C^\times ) \\
\mc R_I = ( X_I ,R_I ,Y_I ,R_I^\vee ,I) & \mc R^I = (X,R_I ,Y
,R_I^\vee ,I)
\end{array}
\end{equation}
We can define parameter functions $q_I$ and $q^I$ on the root
data $\mc R_I$ and $\mc R^I$. Restrict $q$ to a labeling of
$(R_I )_{nr}^\vee$ and use $\eqref{eq:2.18}$ to extend it to
$W (\mc R_I )$ and $W (\mc R^I )$. Then $\mc H (\mc R^I ,q^I )$
is isomorphic to the subalgebra of $\mc H (\mc R ,q)$ generated
by $\mc A$ and $\mc H (W_I ,q)$. With this identification in
mind we call $\mc H (\mc R^I ,q^I )$ a parabolic subalgebra of
$\mc H (\mc R ,q)$.

For any $t \in T^I$ there is a surjective algebra homomorphism
\begin{equation}\label{eq:2.19}
\begin{aligned}
& \phi_t : \mc H (\mc R^I ,q^I ) \to \mc H (\mc R_I ,q_I ) \\
& \phi_t ( \theta_x T_w ) = t(x) \theta_{x_I} T_w
\end{aligned}
\end{equation}
where $x_I$ is the image of $x \in X$ in $X_I$. So given any
representation $\sigma$ of $\mc H (\mc R_I ,q_I )$ we can
construct the $\mc H$-representation
\[
\pi (I,\sigma ,t) := \mr{Ind}_{\mc H (\mc R^I ,q^I )}^{
\mc H (\mc R ,q)} (\sigma \circ \phi_t )
\]
Representations of this form are said to be parabolically induced.

Since $\mc H$ is of finite rank over $Z (\mc H )$ every irreducible
$\mc H$-representation has finite dimension. In particular an
$\mc H$-module is of finite length if and only if it has finite
dimension. Let $\mr{Mod} (\mc H )$ be the category of all
$\mc H$-modules and $\mr{Mod}_{fin} (\mc H )$ the subcategory of
finite length $\mc H$-modules. We denote the Grothendieck group
of $\mr{Mod}_{fin} (\mc H )$ by $G (\mc H )$ and we write
\[
G_{\mh C}(\mc H ) := G (\mc H ) \otimes_{\mh Z} \mh C
\]
Similarly we can define $\mr{Mod} (A) , \mr{Mod}_{fin} (A)
,G (A)$ and $G_{\mh C} (A)$ for any algebra or group $A$.

The center of $\mc H (\mc R, q)$ contains the group algebra of
$Z(W)$, so every irreducible $\mc H$-representation admits a
unique $Z(W)$-character $\chi$. Such representations factor
through the algebra
\[
\mc H (\mc R ,q)_\chi = \mc H \otimes_{Z(W)} \mh C_\chi
\]
The algebra $\mc H$ is endowed with a trace
\[
\tau \Big( \sum_{w \in W} h_w T_w \Big) = h_e
\]
and an involution
\[
\Big( \sum_{w \in W} h_w T_w \Big)^* =
\sum_{w \in W} \overline{h_w} T_{w^{-1}}
\]
Because $q$ takes only positive values, * is conjugate-linear and
antimultiplicative while $\tau$ is positive.

Our affine Hecke algebra is canonically isomorphic to the crossed
product of the Iwahori-Hecke algebra corresponding to $W^\af$,
and the group $\Omega$:
\[
\mc H (\mc R ,q) \cong \mc H (W^\af ,q) \rtimes \Omega
\]
Let $f$ be a facet of the fundamental alcove $A_\es$ and write
\[
\Omega_f := \{ \omega \in \Omega : p_E \omega (f) = f \}
\]
Then $\Omega_f$ acts on $W_f$, so we can define
\[
\mc H (\mc R ,f,q) := \mc H (W_f ,q) \rtimes \Omega_f
\]
By definition $Z(W) \subset \Omega_f$, so
\[
\mh C [Z(W)] \subset Z (\mc H (\mc R ,f,q))
\]

\begin{lem}\label{lem:2.5}
Let $\mh C_\chi$ be a onedimensional $Z(W)$-representation with
character $\chi$.
\[
\mc H (\mc R ,f,q)_\chi :=
\mc H (\mc R ,f,q) \otimes_{Z(W)} \mh C_\chi
\]
is a finite dimensional semisimple algebra.
\end{lem}
\emph{Proof.}
As vector spaces we may identify
\[
\mc H (\mc R ,f,q)_\chi = \mr{Ind}_{\mh C [Z(W)]}^{\mc H (\mc R ,f,q)}
\mh C_\chi = \mc H (W_f ,q) \otimes_{\mh C} \mh C [\Omega_f / Z(W)]
\]
We can extend $|\chi |$ canonically to $X \otimes \mh R$, making
it 1 on $E$. Using this extension we define an involution $*_\chi$
on $\mc H (\mc R, f,q)$ by
\[
(h_w T_w )^{*_\chi} = \overline{h_w} \: |\chi |(2 w (0)) T_{w^{-1}}
\]
The associated bilinear form is
\[
\inp{h}{h'}_\chi = \tau ( h^{*_\chi} \cdot h' )
\]
By construction $\mr{Ind}_{\mh C [Z(W)]}^{\mc H (\mc R ,f,q)}
\mh C_\chi$ is now a unitary representation. This makes
$\mc H (\mc R ,f,q)_\chi$ into a finite dimensional Hilbert
algebra, so in particular it is semisimple. $\qquad \Box$
\\[4mm]

\section{The Schwartz completion}

We introduce the Schwartz completion $\mc S$ of $\mc H$ and discuss
some properties of $\mc S$-modules.

The involution and the trace on $\mc H (\mc R ,q)$ give rise to a
Hermitian inner product
\[
\inp{h}{h'} = \tau (h^* \cdot h') \qquad h,h' \in \mc H (\mc R ,q)
\]
and a norm
\[
\norm{h}_\tau = \sqrt{\inp{h}{h}} = \sqrt{\tau (h^* \cdot h)}
\]
With a basic calculation one can check that
\begin{equation}\label{eq:2.24}
\{ N_w = q(w)^{-1/2} T_w : w \in W \}
\end{equation}
is an orthonormal basis of $\mc H (\mc R ,q)$ for this inner product.
All this gives $\mc H (\mc R ,q)$ the structure of a Hilbert algebra,
in the sense of \cite[A 54]{Dix}.
Let $L^2 (\mc R ,q)$ be its Hilbert space completion, for which
\eqref{eq:2.24} is by definition a basis.
Consider the multiplication map
\begin{align*}
& \lambda (h) : \mc H (\mc R ,q) \to \mc H (\mc R ,q) \\
& \lambda (h) \, h' \;=\; h \cdot h'
\end{align*}
By \cite[Lemma 2.3]{Opd1} this maps extends to a bounded operator on
$L^2 (\mc R ,q)$, whose norm we denote by
\[
\norm{h}_o = \norm{\lambda (h)}_{B (L^2 (\mc R ,q ))}
\]
Thus, $\mc H (\mc R ,q)$ being a *-subalgebra of the $C^*$-algebra
$B (L^2 (\mc R ,q))$ of bounded operators on $L^2 (\mc R ,q)$,
we can consider its closure $C^* (\mc R ,q)$ with respect to the
operator norm topology. By definition this is a separable unital
$C^*$-algebra, called the (reduced) $C^*$-algebra of $\mc H$ or of
$(\mc R ,q)$.

Let $(\pi ,V)$ be an irreducible $\mc H$-representation. We say that
it belongs to the discrete series if the following equivalent
conditions hold:
\begin{itemize}
\item $(\pi ,V)$ is a subrepresentation of the left regular
representation $(\lambda, L^2 (\mc R, q))$
\item all matrix coefficients of $(\pi ,V)$ are in $L^2 (\mc R ,q)$
\end{itemize}

By definition a discrete series representation is unitary, and
it extends continuously to $C^* (\mc R ,q)$. Because this is a
Hilbert algebra, a suitable version of \cite[Proposition 18.4.2]{Dix}
shows that $\pi$ is an isolated point in its spectrum. Moreover, since
$C^* (\mc R ,q)$ is unital its spectrum is compact \cite[Proposition
3.18]{Dix}, so there can be only finitely many inequivalent discrete
series representations.

It is also possible to complete $\mc H (\mc R ,q)$ to a Schwartz
algebra $\mc S$. As a topological vector space
$\mc S$ will consist of rapidly decreasing functions on $W$, with
respect to some length function. For this purpose it is unsatisfactory
that $\ell$ is 0 on the subgroup $Z(W)$, as this can be a large part
of $W$. To overcome this inconvenience, let\\
$L : X \otimes \mh R \to [0,\infty )$ be a function such that
\begin{itemize}
\item $L (X) \subset \mh Z$
\item $L(x+y) = L(x) \quad \forall x \in X \otimes \mh R, y \in E$
\item $L$ induces a norm on
$X \otimes \mh R / E \cong Z(W) \otimes \mh R$
\end{itemize}
Now we define for $w \in W$
\[
\mc N (w) := \ell (w) + L (w(0))
\]
so that
\[
\begin{array}{llll}
\mc N (u \omega ) & = & \mc N (\omega u) \;=\;
\ell (u) + L(\omega (0)) & u \in W^\af , \omega \in \Omega \\
\mc N (w v) & \leq & \mc N (w) + \mc N (v) & w,v \in W
\end{array}
\]
Since $Z(W) \oplus \mh Z R_0$ is of finite index in $X$, the set
$\{ w \in W : \mc N (w) = 0 \}$ is finite. Moreover, because $W$ is
the semidirect product of a finite group and an abelian group, it
is of polynomial growth, and different choices of $L$ lead to
equivalent length functions $\mc N$.
For $n \in \mh N$ we define the norm
\[
p_n \Big( \sum_{w \in W} h_w N_w \Big) :=
\sup_{w \in W} |h_w | (\mc N (w) + 1 )^n
\]
The completion $\mc S = \mc S (\mc R ,q)$ of $\mc H (\mc R ,q)$
with respect to the family of norms $\{ p_n \}_{n \in \mh N}$
is a nuclear Fr\'echet space. It consists of all possible
infinite sums $h = \sum_{w \in W} h_w N_w$ such that
$p_n (h) < \infty \; \forall n \in \mh N$.

\begin{lem}\label{lem:2.10}
\cite[p. 135]{Sol} Let $b = \mr{rk}(X) + 1$. The sum
\[
\sum_{w \in W} \big( \mc N (w) + 1 \big)^{-b}
\]
converges to a limit $C_b$. If $h \in \mc S$ and $n \in \mh N$ then
\[
\sum_{w \in W} |h_w | (\mc N (w) + 1)^n \leq C_b \, p_{n+b}(h)
\]
\end{lem}

\noindent
The norms $p_n$ behave reasonably with respect to multiplication:

\begin{thm}\label{thm:2.7}
\cite[Section 6.2]{Opd1} There exist $C_q > 0 ,\, d \in \mh N$
such that $\forall h,h' \in \mc S (\mc R ,q), n \in \mh N$
\begin{align*}
\norm{h}_o &\leq C_q p_d (h) \\
p_n (h \cdot h') &\leq C_q p_{n+d}(h) p_{n+d}(h')
\end{align*}
In particular $\mc S (\mc R ,q)$ is a unital locally convex
*-algebra, and it is contained in $C^* (\mc R ,q)$.
\end{thm}

The reader is referred to \cite{DeOp} for a study of the algebra
$\mc S$ and its Fourier transform. Notice that as a Fr\'echet space
$\mc S (\mc R ,q)$ does not depend on $q$. The basis
$\{ N_w : w \in W \}$ gives rise to a canonical isomorphism between
$\mc S (\mc R ,q)$ and $\mc S (W)$.

For $\ep \in \mh R$ let $q^\ep$ be the parameter function $q^\ep (w) =
q(w)^\ep$. For every $\ep$ we have the affine Hecke algebra $\mc H
(\mc R ,q^\ep )$ and its Schwartz completion $\mc S (\mc R ,q^\ep )$.
We note that $\mc H (\mc R ,q^0 ) = \mh C [W]$ is the group algebra
of $W$ and that $\mc S (\mc R ,q^0 ) = \mc S (W)$ is the Schwartz
algebra of rapidly decreasing functions on $W$.

The intuitive idea is that these algebras depend continuously on $\ep$.
We will use this in the form of the following rather technical result.

\begin{thm}\label{thm:2.4}
For $\ep \in [-1,1]$ there exists a family of maps
\begin{align*}
& \tilde \sigma_\ep : \mr{Mod}_{fin}(\mc H (\mc R ,q)) \to
\mr{Mod}_{fin}(\mc H (\mc R ,q^\ep )) \\
& \tilde \sigma_\ep (\pi ,V) = (\pi_\ep ,V)
\end{align*}
with the properties
\begin{description}
\item[1)] the map
\[
[-1,1] \to \mr{End}\, V : \ep \to \pi_\ep (N_w)
\]
is analytic for any $w \in W$.
\item[2)] $\tilde \sigma_\ep$ is a bijection if $\ep \neq 0$.
\item[3)] $\tilde \sigma_\ep$ preserves unitarity.
\item[4)] $\tilde \sigma_\ep$ preserves temperedness if $\ep \geq 0$.
\item[5)] $\tilde \sigma_\ep$ preserves the discrete series if $\ep > 0$.
\end{description}
\end{thm}
\emph{Proof.} See \cite[Theorem 5.16 and Lemma 5.17]{Sol}. $\qquad
\Box$

%% file: tempered2.tex
\chapter{Projective resolutions}

In this chapter we will contruct projective resolutions for modules
of an affine Hecke algebra $\mc H$. We do this in a functorial way,
starting from an explicit projective $\mc H$-bimodule resolution of
$\mc H$. This allows us to show that the global dimension
of $\mc H$ equals the rank of the lattice $X$.

It turns out that the same contructions also work over $\mc S$.
However this is by no means automatic. Namely, it is not enough to
have a projective $\mc H$-bimodule resolution, to show that it can
be induced to $\mc S$ we also need a contraction which is bounded in
a suitable sense. The essential part of the proof takes place within the
polysimplicial complex $\Sigma$ associated to the root system $R_0$.
Taking advantage of the abundant symmetry of root systems we
construct a bounded contraction of the corresponding differential
complex. With this contraction we establish a projective bimodule
resolution of $\mc S$. As a consequence we can show that the
cohomological dimension of $\mr{Mod}_{bor} (\mc S )$ also equals the
rank of $X$.

Actually more is true, as Ralf Meyer kindly pointed out to us. The
inclusion of complete, unital, bornological algebras $\mc H \to \mc S$
is isocohomological (in the sense discussed in the appendix).

\section{The bounded contraction of the polysimplicial complex}
\label{sec:2.1}

From the polysimplicial complex $\Sigma$ (cf. page \pageref{p:Sigma})
we construct a differential complex $(C_* (\Sigma ), \partial_* )$.
The vector space in degree $n$ is
\begin{equation}\label{eq:2.27}
C_n (\Sigma ) := \mh C \{ \sigma \in \Sigma : \dim \sigma = n \}
\end{equation}
For every $\sigma$ there is a unique facet $f$ of the fundamental
alcove $A_\es$ such that $\sigma$ is $W^\af$-conjugate to the
closure $\bar f$ of $f$ in $E$. We fix an orientation on all the
facets of $A_\es$ and we decree that the map $w : f \to w f$
preserves orientation. This determines a unique orientation on
every simplex of $\Sigma$. With these conventions we can identify
\begin{equation}\label{eq:2.9}
C_n (\Sigma ) = \bigoplus_{f : \dim f = n} \mh C
\left[ W^\af / W_f \right]
\end{equation}
Clearly $\Sigma$ is the direct product of a number (say $r'$) simplicial
complexes corresponding to the irreducible components of $R_0$. Let
\[
\sigma = \sigma^{(1)} \times \cdots \times \sigma^{(r')}
\]
be a polysimplex of $\Sigma$. Denote the vertices of $\sigma^{(j)}$
by $x_i^{(j)}$, so that we can write
\[
\sigma^{(j)} = \left[ x_0^{(j)},x_1^{(j)},\ldots, x_{d_j}^{(j)} \right]
\]
This defines an orientation on $\sigma^{(j)}$ in the sense that
\[
\left[ x_{\lambda (0)}^{(j)},x_{\lambda (1)}^{(j)},\ldots,
x_{\lambda (d_j )}^{(j)} \right] =
\ep (\lambda ) \left[ x_0^{(j)},x_1^{(j)},\ldots, x_{d_j}^{(j)} \right]
\]
for any $\lambda \in S_{d_j}$. The boundary of $\sigma^{(j)}$ is
defined as
\begin{align*}
& \partial \sigma^{(j)} = \partial \left[ x_0^{(j)},x_1^{(j)},\ldots,
x_{d_j}^{(j)} \right] := \sum_{i=0}^{d_j} (-1)^i \left[ x_0^{(j)},
\ldots, x_{i-1}^{(j)} ,x_{i+1}^{(j)}, \ldots, x_{d_j}^{(j)} \right] \\
& \partial \left[ x_0^{(j)} \right] := 0
\end{align*}
Furthermore we define
\[
\partial_n \sigma = \sum_{j=1}^{r'} (-1)^{d_1 + \cdots + d_{j-1}}
\sigma^{(1)} \times \cdots \times \sigma^{(j-1)} \times \partial
\sigma^{(j)} \times \sigma^{(j+1)} \times \cdots \times \sigma^{(r')}
\]
if dim $\sigma = n > 0$. It is easily verified that this operation
satisfies the usual property $\partial \circ \partial = 0$.
We augment this differential complex by
\[
C_{-1}(\Sigma ) = \mh C
\]
and $\partial_0 [x] = 1$ if $x$ is a vertex of $\Sigma$.
The augmented complex $( C_* (\Sigma ) ,\partial_* )$ computes the
reduced singular homology of the space $E$ underlying $\Sigma$. This
space is contractible, so by the Poincar\'e lemma
\begin{equation}\label{eq:2.2}
H_n ( C_* (\Sigma ) ,\partial_* ) = 0 \qquad \forall n \in \mh Z
\end{equation}
The support of a chain $c = \sum_{\sigma \in \Sigma} c_\sigma \sigma
\in C_* (\Sigma )$ is
\[
\mr{supp} \: c = \bigcup_{\sigma : c_\sigma \neq 0} \sigma
\]
A contraction $\gamma$ of $(C_* (\Sigma ), \partial_* )$ is a
collection of linear maps
\[
\gamma_n : C_n (\Sigma ) \to C_{n+1} (\Sigma ) \qquad n \geq -1
\]
such that
\[
\gamma_{n-1} \partial_n + \partial_{n+1} \gamma_n =
\mr{id}_{C_n (\Sigma )} \qquad \forall n \in \mh Z
\]
The periodic nature of $\Sigma$ allows us to construct a contraction
with good bounds on the coefficients:

\begin{prop}\label{prop:2.3}
There exists a contraction $\gamma$ with the properties
\begin{description}
\item[1)] $\gamma \partial + \partial \gamma =$ id
\item[2)] $\gamma$ is $W_0$-equivariant
\item[3)] supp $\gamma (\sigma ) \subset A (\sigma )$ for every
$\sigma \in \Sigma$
\item[4)] $\gamma (\sigma ) = \sum_{\tau \in \Sigma}
\gamma_{\sigma \tau} \tau$ with $|\gamma_{\sigma \tau}| < M_\gamma$
for some constant $M_\gamma$ depending only on $\gamma$
\end{description}
\end{prop}
\emph{Proof.}
Our construction will be rather similar to that of V. Lafforgue in
\cite[\S 4]{Ska}. First we impose some extra conditions. 2) and 3) force
\begin{description}
\item[5)] if $\sigma \subset C_I^+$
then supp $\gamma (\sigma ) \subset C_I^+$
\end{description}
In view of \eqref{eq:1.was} and since $\partial$ is $W_0$-equivariant,
it suffices to construct $\gamma$
on $C_\es^+$. We will use that the translations $t_x$ with $x \in
\mh Z R_0$ are orientation preserving automorphisms of $\Sigma$. For
$\alpha_i \in F_0$ let $\beta_i$ be the minimal element of
$C^{++}_{F_0 \setminus \{ \alpha_i \}} \cap \mh Z R_0$. Note that
$\beta_i$ is an integral multiple of a vertex of $A_\es$. We could
also pick a fundamental weight instead of $\beta_i$, but in that
case we would have keep track of the orientations.
Consider the halfopen parallelogram
\[
P_\es = \Big\{ \sum_{i=1}^r y_i \beta_i : y_i \in [0,1) \Big\}
\]
Let $\tau$ be any polysimplex whose interior is contained in $P_\es$.
Our contraction will also satisfy
\begin{description}
\item[6)] $\gamma ( t_{(m+1) \beta_i} (\tau )) = \gamma (t_{m \beta_i}
(\tau ) ) + t_{m \beta_i} \gamma (t_{\beta_i} (\tau ) - \tau )$
\end{description}
for $m \geq 0$. Suppose that $\beta = \sum_{i=1}^k n_i \beta_i$ with
$n_i \in \mh N$. Then we decree
\begin{description}
\item[7)] $\gamma ( t_\beta (\tau )) = \gamma ( t_{n_k \beta_k} (\tau ) )
+ t_{n_k \beta_k} \gamma ( t_{\beta - n_k \beta_k} (\tau ) - \tau )$
\end{description}
Here we use the ordering on the set $F_0$ of simple roots. The idea
underlying 6) and 7) is that we want to make $\gamma$ equivariant with
respect to certain translations.

Now we really start constructing $\gamma$. In degree $-1$ we put
\[
\gamma_{-1}(1) = [0]
\]
Suppose that $\gamma_m$ has already been defined for $m < n$, satisfying
conditions 1) - 7). Let $\sigma$ be any $n$-dimensional polysimplex
whose interior is contained in
\[
P_1 := P_\es \cup t_{\beta_1} P_\es \cup \cdots \cup t_{\beta_r} P_\es
\]
By 1) we have
\[
\partial (\sigma - \gamma \partial (\sigma )) =
(\mr{id} - \partial \gamma ) (\partial \sigma ) =
\gamma \partial (\partial \sigma ) = 0
\]
Together with \eqref{eq:2.2} this implies that the equation
\[
\partial \gamma (\sigma ) = \sigma - \gamma \partial (\sigma )
\]
has a solution $\gamma (\sigma ) \in C_{n+1}(\Sigma )$.
By 3) and 5) we have
\[
\mr{supp} \: (\sigma - \gamma \partial (\sigma )) \subset
A(\sigma ) \cap C_I^+ \quad \mr{if} \quad \sigma \subset C_I^+
\]
Since $A(\sigma ) \cap C_I^+$ is convex, we
can pick $\gamma (\sigma )$ with support in this set. We do this for
any $n$-dimensional $\sigma \in \Sigma$ whose interior is contained in
$P_1$. Now 6) and 7) determine $\gamma_n$ uniquely on $C_\es^+$.

We will show that the other required properties follow from this
construction. Write $\beta' = \sum_{i=1}^{k-1} n_i \beta_i$ and
$\beta'' = \sum_{i=1}^{k-1} n'_i \beta_i$ for some $n'_i \in \mh N$.
By 7) we have
\begin{equation}\label{eq:2.11}
\gamma t_{n_k \beta_k} \big( t_{\beta'} (\tau ) -
t_{\beta''} (\tau ) \big) = t_{n_k \beta_k} \gamma \big( t_{\beta'}
(\tau ) - t_{\beta''} (\tau ) \big)
\end{equation}
We claim that the following stronger version of 7) holds
\begin{description}
\item[7')] $\gamma t_{n_k \beta_k} \big( t_{\beta - n_k \beta_k}
(\sigma ) - \sigma \big) = t_{n_k \beta_k} \gamma \big( t_{\beta - n_k
\beta_k} (\sigma ) - \sigma \big) \qquad \forall \sigma \subset C^+_\es$
\end{description}
Indeed, write $\sigma = t_x \tau$ with $\tau$ as in 7) and
$x = \sum_{j=1}^r m_j \beta_j$. Then by a repeated application of
\eqref{eq:2.11} the left hand side of 7') becomes
\begin{align*}
\gamma t_{n_k \beta_k} \big( t_{\beta'} (t_x \tau ) - t_x \tau \big) \,
& =\: t_{(n_k + m_k) \beta_k + m_{k+1} \beta_{k+1} + \cdots + m_r
\beta_r} \gamma (t_{\beta'} - \mr{id} ) t_{m_1 \beta_1 + \cdots + m_{k-1}
\beta_{k-1}} (\tau ) \\
& =\: t_{n_k \beta_k} \gamma t_x (t_{\beta'} (\tau ) - \tau ) \\
& =\: t_{n_k \beta_k} \gamma (t_{\beta'} (\sigma ) - \sigma )
\end{align*}
It follows easily from 6) that
\begin{equation}\label{eq:2.12}
\gamma t_{m \beta_i} \big( t_{m' \beta_i} (\tau ) - t_{m'' \beta_i}
(\tau ) \big) = t_{m \beta_i} \gamma \big( t_{m' \beta_i} (\tau ) -
t_{m'' \beta_i} (\tau ) \big) \qquad \forall m, m', m'' \in \mh N
\end{equation}
There also is a stronger version of 6) :
\begin{description}
\item[6')] $\gamma t_{m \beta_i} \big( t_{\beta_i} (\sigma ) - \sigma
\big) = t_{m \beta_i} \gamma \big( t_{\beta_i} (\sigma ) - \sigma \big)
\qquad \forall \sigma \subset C^+_\es$
\end{description}
Indeed, in the above notation and by 7') and \eqref{eq:2.12}
the left hand side equals
\[
\begin{array}{ll}
\gamma t_{m \beta_i} \big( t_{\beta_i + x} (\tau ) - t_x (\tau ) \big)
& = \\
t_{m_{i+1} \beta_{i+1} + \cdots + m_r \beta_r} \gamma
t_{(m+m_i ) \beta_i} (t_{\beta_i} - \mr{id}) t_{m_1 \beta_1 + \cdots +
m_{i-1} \beta_{i-1}} (\tau ) & = \\
t_{m_{i+1} \beta_{i+1} + \cdots + m_r \beta_r} \Big( \gamma
t_{(m+m_i ) \beta_i} \big( t_{\beta_i} (\tau ) - \tau \big) & \\
\quad +\; t_{(m+m_i ) \beta_i} (t_{\beta_i} - \mr{id}) \gamma \big( t_{m_1
\beta_1 + \cdots + m_{i-1} \beta_{i-1}} (\tau ) - \tau \big) \Big) & = \\
t_{m \beta_i + m_{i+1} \beta_{i+1} + \cdots + m_r \beta_r} \Big(
\gamma t_{m_i \beta_i} \big( t_{\beta_i} (\tau ) - \tau \big) & \\
\quad +\; t_{m_i \beta_i} (t_{\beta_i} - \mr{id}) \gamma \big( t_{m_1
\beta_1 + \cdots + m_{i-1} \beta_{i-1}} (\tau ) - \tau \big) \Big) & = \\
t_{m \beta_i + m_{i+1} \beta_{i+1} + \cdots + m_r \beta_r}
\gamma (t_{\beta_i} - \mr{id}) t_{m_1 \beta_1 + \cdots + m_i \beta_i}
(\tau ) & = \\
t_{m \beta_i} \gamma (t_{\beta_i} - \mr{id}) t_x (\tau ) & =\;
t_{m \beta_i} \gamma \big( t_{\beta_i} (\sigma ) - \sigma \big)
\end{array}
\]
Now we can see that the relations 6) and 7) are compatible with 1).
Assume that 1) holds for $t_{m \beta_i} (\tau )$. Then by 6')
\[
\begin{array}{ll}
( \partial_{n+1} \gamma_n + \gamma_{n-1} \partial_n ) \big( t_{(m+1)
\beta_i} (\tau ) \big) & = \\
\partial_{n+1} \gamma_n (t_{m \beta_i} \tau ) + \partial_{n+1}
t_{m \beta_i} \gamma_n \big( t_{\beta_i} (\tau ) - \tau \big) +
\gamma_{n-1} t_{(m+1) \beta_i} \partial_n (\tau ) & = \\
\partial_{n+1} \gamma_n (t_{m \beta_i} \tau ) + t_{m \beta_i}
\partial_{n+1} \gamma_n \big( t_{\beta_i} (\tau ) - \tau \big) & \\
\quad +\; \gamma_{n-1} t_{m \beta_i} \partial_n (\tau ) + t_{m \beta_i }
\gamma_{n-1} \big( t_{\beta_i} (\partial_n \tau ) - \partial_n \tau \big)
& = \\
t_{m \beta_i} (\tau ) + t_{m \beta_i} \big( t_{\beta_i} (\tau ) - \tau
\big) & =\; t_{(m+1) \beta_i} (\tau )
\end{array}
\]
Similarly, suppose that $t_{n_k \beta_k} (\sigma )$ and
$t_{\beta - n_k \beta_k} (\sigma )$ both satisfy 1).
It follows from 7') that
\[
\begin{array}{ll}
(\partial_{n+1} \gamma_n + \gamma_{n-1} \partial_n )
\big( t_\beta (\sigma ) \big) & = \\
\partial_{n+1} \gamma_n (t_{n_k \beta_k} (\sigma )) + \partial_{n+1}
(t_{n_k \beta_k} \gamma_n (t_{\beta - n_k \beta_k} (\sigma ) - \sigma )))
+ \gamma_{n-1} (t_\beta \partial_n (\sigma )) & = \\
\partial_{n+1} \gamma_n (t_{n_k \beta_k} (\sigma )) + t_{n_k \beta_k}
\partial_{n+1} \gamma_n (t_{\beta - n_k \beta_k} (\sigma ) - \sigma ))
& \\
\quad + \; \gamma_{n-1} (t_{n_k \beta_k} \partial_n (\sigma )) +
t_{n_k \beta_k} \gamma_{n-1} (t_{\beta - n_k \beta_k} (\partial_n \sigma )
- \partial_n (\sigma ))) & = \\
t_{n_k \beta_k} (\sigma ) + t_{n_k \beta_k} (t_{\beta - n_k \beta_k}
(\sigma ) - \sigma ) ) & =\; t_\beta (\sigma )
\end{array}
\]
Thus we can construct $\gamma$ respecting all conditions, except
possibly 3) and 4). The parallelogram $P_2 = 2 \overline P_\es$ consists
of finitely many polysimplices, so there is a real number $M$ such that
\[
\gamma (\tau ) = \sum_{\sigma} \gamma_{\tau \sigma} \sigma
\quad \mr{with} \quad |\gamma_{\tau \sigma}| < M
\]
for all polysimplices $\tau \subset P_2$. Let us examine the size of
the coefficients of $\gamma (t_{(m+1) \beta_i} (\sigma ))$ for $\tau$
with interior in $P_\es$. By induction to $m$ we may suppose that
\begin{equation}\label{eq:2.3}
\gamma (t_{m \beta_i} (\tau )) = \sum_{\sigma} \lambda_\sigma^m \sigma
\quad \mr{with} \quad |\lambda_\sigma^m |
\left\{ \begin{array}{cccccl}
 = & 0 & \mr{if} & \sigma & \not\subset & A(t_{m \beta_i} (\tau )) \\
 < & M & \mr{if} & \sigma & \subset & P_2 \\
 < & M & \mr{if} & \sigma & \not\subset & A(t_{(m-1) \beta_i} (\tau )) \\
 < & 3M & \mr{if} & \sigma & \subset & A(t_{(m-1) \beta_i} (\tau ))
 \end{array} \right.
\end{equation}
By construction we have
\[
t_{m \beta_i} \gamma (t_{\beta_i} (\tau ) - \tau ) = \sum_{\sigma}
\lambda'_\sigma \sigma \quad \mr{with} \quad |\lambda'_\sigma |
\left\{ \begin{array}{cccccl}
= & 0 & \mr{if} & \sigma & \not\subset & A (t_{(m+1)\beta_i} (\tau )) \\
= & 0 & \mr{if} & \sigma & \not\subset & t_{m \beta_i} C_\es^+ \\
< & M & \mr{if} & \sigma & \not\subset & A( t_{m \beta_i} (\tau )) \\
< & 2M & \mr{if} & \sigma & \subset & A( t_{(m+1) \beta_i} (\tau ))
 \end{array} \right.
\]
With 6) this implies that \eqref{eq:2.3} also holds with $m+1$ instead
of $m$.

Let $\beta$ be as above. By induction to $k$ we may assume that
\begin{equation}\label{eq:2.4}
t_{n_k \beta_k} \gamma (t_{\beta - n_k \beta_k} (\tau ) - \tau ) =
\sum_\sigma \mu^k_\sigma \sigma \quad \mr{with} \quad |\mu^k_\sigma |
\left\{ \begin{array}{cccccl}
= & 0 & \mr{if} & \sigma & \not\subset & A (t_\beta (\sigma )) \\
= & 0 & \mr{if} & \sigma & \not\subset & t_{n_k \beta_k} C_\es^+ \\
< & M & \mr{if} & \sigma & \not\subset & t_{\beta'} (\sigma ) \\
< & 2M & \mr{if} & \sigma & \subset & A (t_{n_k \beta_k} (\sigma )) \\
< & 3M & \mr{if} & \sigma & \subset & A (t_\beta (\sigma ))
 \end{array} \right.
\end{equation}
where $\beta' = \beta - \beta_i$ with $i$ minimal for $n_i > 0$. In view
of 7) the above implies that
\[
\gamma (t_\beta (\tau )) = \sum_\sigma \mu'_\sigma \sigma \quad \mr{with}
\quad |\mu'_\sigma | \left\{ \begin{array}{cccccl}
= & 0 & \mr{if} & \sigma & \not\subset & A (t_\beta (\sigma )) \\
< & M & \mr{if} & \sigma & \subset & P_2 \\
< & M & \mr{if} & \sigma & \not\subset & A (t_{\beta'} (\sigma )) \\
< & 3M & \mr{if} & \sigma & \subset & A (t_\beta (\sigma ))
 \end{array} \right.
\]
This in turn implies \eqref{eq:2.4} with $k+1$ instead of $k$. Hence
condition 4) is fulfilled, with $M_\gamma = 3M. \qquad \Box$
\\[3mm]
\textbf{Example.}\\
In the case $R_0 = B_2$ we have $\beta_1 = (1,0)$ and $\beta_2 = (1,1)$.
We drew the sets $P_\es ,\, P_1$ and $P_2$ below.
If $x$ is a vertex of $\Sigma$ then $\gamma [x]$ is a path from 0 to $x$,
along the following lines:
\\[2mm]

\includegraphics[width=115mm,height=49mm]{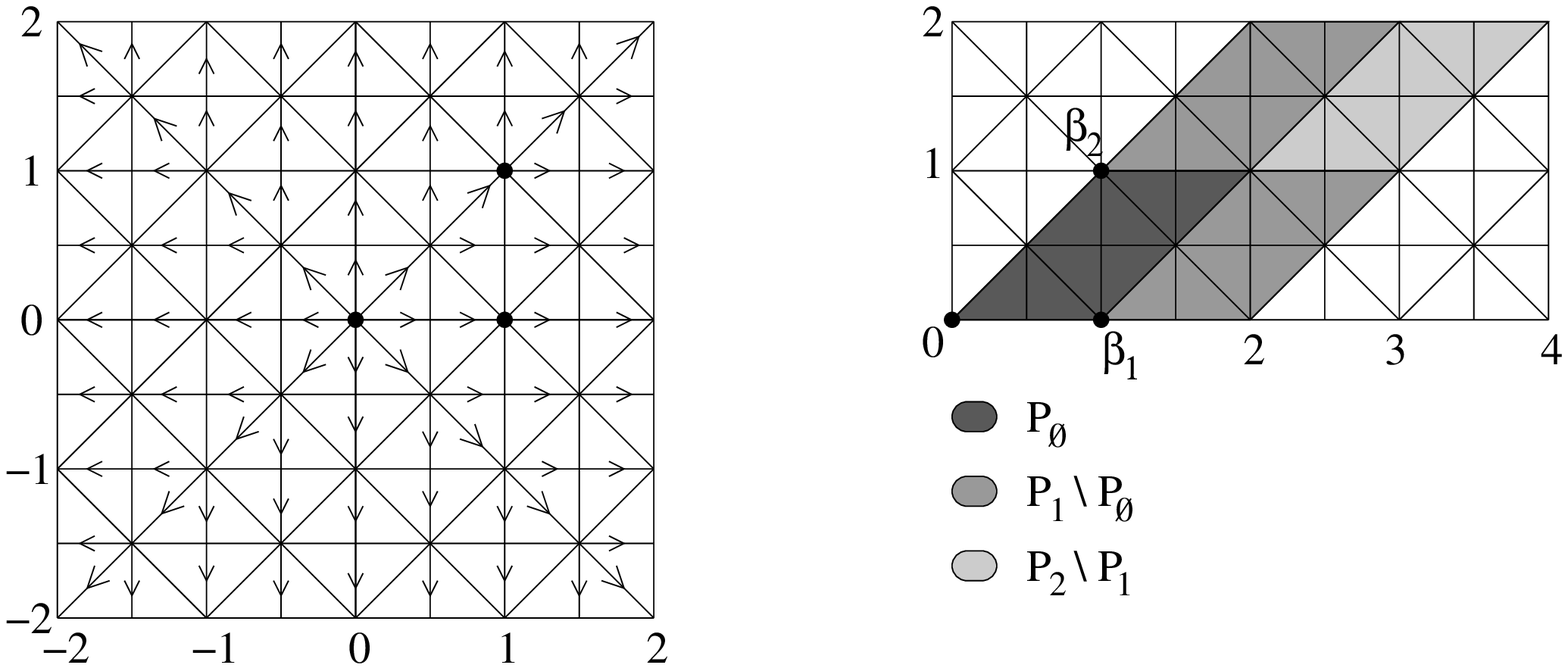}
\\[2mm]
We define
\[
\begin{array}{llcll}
\gamma \, [(1/2,0),(1/2,1/2)] & = & A_\es &
 = & [(0,0),(1/2,0),(1/2,1/2)] \\
\gamma \, [(1,1/2),(1,1)] & = & t_{(1/2,1/2)} A_\es &
 = & [(1/2,1/2),(1,1/2),(1,1)] \\
\gamma \, [(3/2,1),(3/2,3/2)] & = & t_{(1,1)} A_\es &
 = & [(1,1),(3/2,1),(3/2,3/2)] \\
\gamma \, [(3/2,0),(3/2,1/2)] & = & \multicolumn{3}{l}{

\begin{picture}(4,1)(-1,0)
\put(0,0){\line(1,0){1.5}}
\put(0.5,0.5){\line(1,0){1}}
\put(0,0){\line(1,1){0.5}}
\put(0.5,0.5){\line(1,-1){0.5}}
\put(1,0){\line(1,1){0.5}}
\put(0.5,0){\line(0,1){0.5}}
\put(1,0){\line(0,1){0.5}}
\put(1.5,0){\line(0,1){0.5}}
\put(-0.2,-0.2){\makebox{\scriptsize 0}}
\put(-0.4,0.3){\makebox{\scriptsize{1/2}}}
\put(1.4,-0.3){\makebox{\scriptsize{3/2}}}
\end{picture}}

\end{array}
\]
According to 6)
\begin{align*}
& \gamma \, [(5/2,0),(5/2,1/2)] = \\
& \gamma \, [(3/2,0),(3/2,1/2)] + t_{(1,0)} \gamma \,
\Big( [(3/2,0),(3/2,1/2)] - [(1/2,0),(1/2,1/2)] \Big) =
\end{align*}

\begin{picture}(13,1)(-0.4,-0.6)
\put(0,0){\line(1,0){1.5}}
\put(0.5,0.5){\line(1,0){1}}
\put(0,0){\line(1,1){0.5}}
\put(0.5,0.5){\line(1,-1){0.5}}
\put(1,0){\line(1,1){0.5}}
\put(0.5,0){\line(0,1){0.5}}
\put(1,0){\line(0,1){0.5}}
\put(1.5,0){\line(0,1){0.5}}
\put(-0.2,-0.2){\makebox{\scriptsize 0}}
\put(-0.4,0.3){\makebox{\scriptsize{1/2}}}
\put(1.4,-0.3){\makebox{\scriptsize{3/2}}}

\put(1.9,0.1){\makebox{$+ \;\, t_{(1,0)} \, \bigg( $}}

\put(4,0){\line(1,0){1.5}}
\put(4.5,0.5){\line(1,0){1}}
\put(4,0){\line(1,1){0.5}}
\put(4.5,0.5){\line(1,-1){0.5}}
\put(5,0){\line(1,1){0.5}}
\put(4.5,0){\line(0,1){0.5}}
\put(5,0){\line(0,1){0.5}}
\put(5.5,0){\line(0,1){0.5}}
\put(3.8,-0.2){\makebox{\scriptsize 0}}
\put(3.6,0.3){\makebox{\scriptsize{1/2}}}
\put(5.4,-0.3){\makebox{\scriptsize{3/2}}}

\put(5.9,0.2){\line(1,0){0.25}}

\put(6.7,0){\line(1,0){0.5}}
\put(6.7,0){\line(1,1){0.5}}
\put(7.2,0){\line(0,1){0.5}}
\put(6.5,-0.2){\makebox{\scriptsize 0}}
\put(6.3,0.3){\makebox{\scriptsize{1/2}}}
\put(7.1,-0.3){\makebox{\scriptsize{1/2}}}

\put(7.7,0.1){\makebox{$\bigg) \: =$}}

\put(9,0){\line(1,0){2.5}}
\put(9.5,0.5){\line(1,0){2}}
\put(9,0){\line(1,1){0.5}}
\put(9.5,0.5){\line(1,-1){0.5}}
\put(10,0){\line(1,1){0.5}}
\put(10.5,0.5){\line(1,-1){0.5}}
\put(11,0){\line(1,1){0.5}}
\put(9.5,0){\line(0,1){0.5}}
\put(10,0){\line(0,1){0.5}}
\put(10.5,0){\line(0,1){0.5}}
\put(11,0){\line(0,1){0.5}}
\put(11.5,0){\line(0,1){0,5}}
\put(8.8,-0.2){\makebox{\scriptsize 0}}
\put(8.6,0.3){\makebox{\scriptsize{1/2}}}
\put(11.4,-0.3){\makebox{\scriptsize{5/2}}}

\end{picture}
\\
Condition 7) says that
\begin{align*}
& \gamma \, [(7/2,1),(7/2,3/2)] = \\
& \gamma \, [(3/2,1),(3/2,3/2)] + t_{(1,1)} \gamma \,
\Big( [(5/2,0),(5/2,1/2)] - [(1/2,0),(1/2,1/2)] \Big) =
\end{align*}

\begin{picture}(13,1.2)(-0.4,-0.8)
\put(0,0){\line(1,0){0.5}}
\put(0,0){\line(1,1){0.5}}
\put(0.5,0){\line(0,1){0.5}}
\put(-0.25,0.4){\makebox{\scriptsize 3/2 }}
\put(0.35,-0.3){\makebox{\scriptsize 3/2 }}
\put(-0.25,0){\makebox{\scriptsize 1 }}
\put(-0.05,-0.3){\makebox{\scriptsize 1 }}

\put(1,0.1){\makebox{$+ \;\; t_{(1,1)} \; \bigg( $}}

\put(3.5,0){\line(1,0){2.5}}
\put(4,0.5){\line(1,0){2}}
\put(3.5,0){\line(1,1){0.5}}
\put(4,0.5){\line(1,-1){0.5}}
\put(4.5,0){\line(1,1){0.5}}
\put(5,0.5){\line(1,-1){0.5}}
\put(5.5,0){\line(1,1){0.5}}
\put(4,0){\line(0,1){0.5}}
\put(4.5,0){\line(0,1){0.5}}
\put(5,0){\line(0,1){0.5}}
\put(5.5,0){\line(0,1){0.5}}
\put(6,0){\line(0,1){0,5}}
\put(3.3,-0.2){\makebox{\scriptsize 0}}
\put(3.1,0.3){\makebox{\scriptsize{1/2}}}
\put(5.9,-0.3){\makebox{\scriptsize{5/2}}}

\put(6.4,0.2){\line(1,0){0.25}}

\put(7.3,0){\line(1,0){0.5}}
\put(7.3,0){\line(1,1){0.5}}
\put(7.8,0){\line(0,1){0.5}}
\put(7.1,-0.2){\makebox{\scriptsize 0}}
\put(6.9,0.3){\makebox{\scriptsize{1/2}}}
\put(7.65,-0.3){\makebox{\scriptsize{1/2}}}

\put(8.27,0.1){\makebox{$\bigg) \: =$}}

\put(9.7,0){\line(1,0){2.5}}
\put(10.2,0.5){\line(1,0){2}}
\put(9.7,0){\line(1,1){0.5}}
\put(10.2,0.5){\line(1,-1){0.5}}
\put(10.7,0){\line(1,1){0.5}}
\put(11.2,0.5){\line(1,-1){0.5}}
\put(11.7,0){\line(1,1){0.5}}
\put(10.2,0){\line(0,1){0.5}}
\put(10.7,0){\line(0,1){0.5}}
\put(11.2,0){\line(0,1){0.5}}
\put(11.7,0){\line(0,1){0.5}}
\put(12.2,0){\line(0,1){0,5}}
\put(9.45,0){\makebox{\scriptsize 1 }}
\put(9.65,-0.3){\makebox{\scriptsize 1 }}
\put(9.4,0.4){\makebox{\scriptsize 3/2 }}
\put(12.1,-0.3){\makebox{\scriptsize 7/2}}
\end{picture}

\section{Projective resolutions for affine Hecke algebras}

For $(\pi ,V) \in \mr{Mod}(\mc H )$ and
$n \in \mh N$ we consider the $\mc H$-module
\[
P_n (V) := \bigoplus_{f : \dim f = n} \mc H \otimes_{\mc H (W_f ,q)
\otimes \mh C [Z(W)]} V \otimes_{\mh C} \mh C \{ f \} =
\bigoplus_{f : \dim f = n} \mc H \otimes_{\mc H (W_f ,q) \otimes \mh C [Z(W)]} V
\]
where the sum runs over facets of $A_\es$. Recall that we already
fixed an (arbitrary) orientation of all these facets. Write
\[
\partial \big( \,\overline f \, \big) =
\sum_{f'} [f : f'] \, \overline{f'}
\]
and define $\mc H$-module homomorphisms
\begin{align}\label{eq:2.14}
& d_n : P_n (V) \to P_{n-1}(V) \\
& d_n \big( h \otimes_{\mc H (W_f ,q) \otimes \mh C [Z(W)]} v \otimes_{\mh C} f
\big) = \sum_{f' : \dim f' = n-1} h \otimes_{\mc H (W_{f'} ,q) \rtimes
Z(W)} v \otimes_{\mh C} [f : f'] f' \nonumber
\end{align}
Furthermore we define
\begin{align}
& d_0 : P_0 (V) \to V \\
& d_0 \big(  h \otimes_{\mc H (W_x ,q) \otimes \mh C [Z(W)]} v
\otimes_{\mh C} x \big) = \pi (h) v \nonumber
\end{align}
if $x$ is a vertex of $A_\es$. Now $\big( P_* (V) ,d_* \big)$
is an augmented differential complex because $\partial \circ \partial =
0$. The group $\Omega$ acts naturally on this complex by
\[
\omega (h \otimes_{\mc H (W_f ,q) \otimes \mh C [Z(W)]} v \otimes f)
= h T_\omega^{-1} \otimes_{\mc H (W_{\omega (f)},q) \otimes \mh C [Z(W)]}
\pi (T_\omega ) v \otimes \omega (f )
\]
This action commutes with the $\mc H$-action and with the differentials
$d_n$, so $\big( P_* (V)^\Omega ,d_* \big)$ is again an augmented
differential complex. Note that $P_n (V)$ and $P_n (V)^\Omega$ are
finitely generated $\mc H$-modules if $V$ has finite dimension.

\begin{thm}\label{thm:2.6}
Consider $\mc H$ as a $\mc H$-bimodule.
\begin{equation}\label{eq:2.25}
0 \longleftarrow \mc H \xleftarrow{\; d_0 \;}  P_0 (\mc H )^\Omega
\xleftarrow{\; d_1 \;} P_1 (\mc H )^\Omega \longleftarrow \cdots
\xleftarrow{\; d_r \;} P_r (\mc H )^\Omega \longleftarrow 0
\end{equation}
is a resolution of $\mc H$ by $\mc H \otimes \mc H^{op}$-modules. Every
$P_n (\mc H )^\Omega$ is projective as a left and as a right $\mc H$-module.
Moreover if $\mc R$ is semisimple then $P_n (\mc H )^\Omega$ is projective
as a $\mc H \otimes \mc H^{op}$-module.
\end{thm}
\emph{Proof.}
This result stems from joint work of Mark Reeder and the first author,
see \cite[Proposition 8.1]{Opd2}. The proof is based on constructions
of Kato \cite{Kat1}.

First we consider the case $\Omega = Z(W) = \{e\} ,\, W =
W^\af$. There is a linear bijection
\begin{equation}\label{eq:2.5}
\begin{aligned}
& \phi : \mh C [W] \otimes_{\mh C} \mc H \to
\mc H \otimes_{\mh C} \mc H\\
& \phi (w \otimes h') = T_w \otimes T_w^{-1} h'
\end{aligned}
\end{equation}
For $s_i \in S_{\mr{aff}}$ we write $q_i = q (s_i )$ and
\begin{equation}\label{eq:2.6}
\begin{array}{cclcl}
L_i & := & \mr{span}\{ h T_{s_i} \otimes T_{s_i}^{-1} h'
- h \otimes h' : h,h' \in \mc H \}
& \subset & \mc H \otimes_{\mh C} \mc H \\
\mh C [W]_i & := & \big\{ \sum_{w \in W}
x_w w : x_{w s_i} = -x_w \; \forall w \in W \big\}
& \subset & \mh C [W]
\end{array}
\end{equation}
This $L_i$ is interesting because
\[
\mc H \otimes_{\mc H (W_f ,q)} \mc H = \big( \mc H \otimes_{\mh C}
\mc H \big) \Big/ \sum_{s_i \in W_f} L_i
\]
Let $w \in W$ be such that $\ell (w s_i) > \ell (w)$.
For any $h' \in \mc H$ we have
\[
\phi ((w s_i - w) \otimes h') = T_{w s_i} \otimes T_{w s_i}^{-1} h'
- T_w \otimes T_w^{-1} h' = T_w T_{s_i} \otimes T_{s_i}^{-1}
T_w^{-1} h' - T_w \otimes T_w^{-1} h' \in L_i
\]
so $\phi (\mh C [W]_i \otimes \mc H) \subset L_i$. On the other hand,
$L_i$ is spanned by elements as in \eqref{eq:2.6} with
$h = T_w$ or $h = T_{w s_i}$.
\[
\begin{array}{ll}
\phi^{-1}(T_{w s_i} T_{s_i} \otimes T_{s_i}^{-1} h' -
T_{w s_i} \otimes h') & = \\
\phi^{-1}(q_i T_w + (q_i -1) T_{w s_i} \otimes
T_{s_i}^{-1} h') - w s_i \otimes T_{w s_i} h' & = \\
q_i w \otimes T_w T_{s_i}^{-1} h' + (q_i -1) w s_i
\otimes T_{w s_i} T_{s_i}^{-1} h' -
w s_i \otimes T_{w s_i} h' & = \\
q_i (w - w s_i ) \otimes T_w T_{s_i}^{-1} h' +
w s_i \otimes \big( q_i T_w T_{s_i}^{-1} + (q_i -1) T_{w s_i}
T_{s_i}^{-1} - T_{w s_i} \big) h' & = \\
(w - w s_i ) \otimes T_w q_i T_{s_i}^{-1} h' +
w s_i \otimes \big( T_w (T_{s_i} + 1 - q_i ) + (q_i -1) T_w
- T_w T_{s_i} \big) h' & = \\
(w - w s_i ) \otimes T_w (T_{s_i}+1 - q_i ) h'
\qquad \in \qquad \mh C [W]_i \otimes \mc H
\end{array}
\]
We conclude that $\phi^{-1}(L_i ) = \mh C [W]_i \otimes \mc H$.
Now we bring the linear bijections
\begin{equation}
\mh C [W] \Big/ \sum_{s_i \in W_f} \mh C [W]_i \to \mh C [W / W_f] :
w \to w W_f
\end{equation}
into play. Under these identifications our differential complex
becomes
\[
0 \leftarrow \mc H \leftarrow \cdots \leftarrow \!
\bigoplus_{f : \dim f = n} \!\!\! \mh C [W / W_f ] \otimes \mc H
\otimes \mh C \{ f \} \leftarrow \cdots \leftarrow \mh C [W]
\otimes \mc H \otimes \mh C \{ A_\es \} \leftarrow 0
\]
But this is just the complex $\big( C_* (\Sigma ), \partial_* \big)$
tensored with $\mc H$, so by \eqref{eq:2.2} its homology vanishes.
This shows that indeed we have a resolution in the special case
$\Omega = \{ e \}$.

Now the general case. Since the action of $\Omega$ on $A_\es$
factors through the finite group $\Omega / Z(W)$
we can construct a Reynolds operator
\[
R_\Omega := [\Omega : Z(W) ]^{-1} \sum_{\omega \in \Omega / Z(W)}
\omega \quad \in \quad
\mr{End}_{\mc H \otimes \mc H^{op}} \big( P_n (\mc H) \big)
\]
Since this is an idempotent
\begin{equation}\label{eq:2.7}
P_n (\mc H )^\Omega = R_\Omega \cdot P_n (\mc H)
\end{equation}
is a direct summand of $P_n (\mc H)$.
We generalize \eqref{eq:2.5} to a bijection
\begin{equation}
\begin{aligned}
& \phi : \mh C [W / Z(W)] \otimes_{\mh C} \mc H \to
\mc H \otimes_{\mh C [Z(W)]} \mc H \\
& \phi (w \otimes h') = T_w \otimes T_w^{-1} h'
\end{aligned}
\end{equation}
Just as above this leads to bijections
\[
\bigoplus_{f : \dim f = n} \mh C [W/(W_f \times Z(W))]
\otimes \mc H \otimes \mh C \{ f \} \to P_n (\mc H)
\]
Since both sides are free $\Omega / Z(W)$-modules we also get a
linear bijection
\begin{equation}\label{eq:2.8}
\begin{aligned}
& \bigoplus_{f : \dim f = n} \mh C \big[ W^\af / W_f \big]
\otimes \mc H \otimes \mh C \{ f \} \to P_n (\mc H)^\Omega \\
& w \otimes h' \otimes f \to R_\Omega \big( T_w \otimes_{
\mc H (W_f ,q) \otimes \mh C [Z(W)]} T_w^{-1} h' \otimes f \big)
\end{aligned}
\end{equation}
Now the same argument as in the special case shows that the
modules $P_n (\mc H)^\Omega$ form a resolution of $\mc H$.

For any facet $f \: \mc H$ is a free $\mc H (W_f ,q) \otimes
\mh C [Z(W)]$-module, both from the left and from the right.
Therefore every $P_n (\mc H )$ is a projective $\mc H$-module,
from the left and from the right.

For $\mc R$ semisimple $P_n (\mc H )$ is a direct sum
of $\mc H \otimes \mc H^{op}$-modules of the form \\
$\mc H \otimes_{\mc H (W_f ,q)} \mc H$. For every irreducible
representation $V_i$ of $\mc H (W_f ,q)$ we pick an idempotent
$e_i \in \mc H (W_f ,q)$ which acts as a rank one projection
on $V_i$ and as 0 on all other irreducible representations.
Consider the element $e_f = \sum_i e_i \otimes e_i \in
\mc H \otimes \mc H^{op}$. From
\begin{equation}\label{eq:2.ef}
\mc H \otimes_{\mc H (W_f ,q)} \mc H \cong
\big( \mc H \otimes_{\mh C} \mc H^{op} \big) e_f
\end{equation}
we see that $P_n (\mc H )$ is a projective $\mc H$-bimodule.
By \eqref{eq:2.7} $P_n (\mc H )^\Omega$ is projective in the
same senses as $P_n (\mc H ) \qquad \Box$
\\[3mm]

\begin{cor}\label{cor:2.11}
\begin{description}
\item[a)] Let $V$ be any $\mc H$-module.
\[
0 \longleftarrow V \xleftarrow{\; d_0 \;} P_0 (V )^\Omega
\xleftarrow{\; d_1 \;} P_1 (V )^\Omega \longleftarrow \cdots
\xleftarrow{\; d_r \;} P_r (V )^\Omega \longleftarrow 0
\]
is a resolution of $V$. It is bornological if $V$ is.
\item[b)] If $V$ admits a $Z(W)$-character $\chi$ then every
$P_n (V)^\Omega$ is a projective $\mc H (\mc R ,q )_\chi$-module.
\item[c)] The cohomological dimensions of $\mr{Mod} \big(
\mc H (\mc R ,q )_\chi \big)$ and  $\mr{Mod}_{bor} \big(
\mc H (\mc R ,q )_\chi \big)$ equal $r = \mr{rk}(R_0 )$.
\end{description}
\end{cor}
\emph{Proof.} a)
Apply $\otimes_{\mc H} V$ to \eqref{eq:2.25}. The resulting
differential complex is exact because $\mc H$ and
$P_n (\mc H )^\Omega$ are projective right $\mc H$-modules.
For $V \in \mr{Mod}_{bor}(\mc H )$ this clearly gives a
bornological differential complex. It is split exact because
every contraction of $P_* (\mc H )^\Omega$ yields a bounded
splitting of $P_* (V)^\Omega$. \\
b) From
\begin{equation}\label{eq:2.ind}
\mc H \otimes_{\mc H (W_f,q) \otimes \mh C [Z(W)]} V \cong
\mc H (\mc R ,q )_\chi \otimes_{\mc H (W_f ,q)} V \cong
\mr{Ind}_{\mc H (W_f ,q )}^{\mc H (\mc R ,q )_\chi } V
\end{equation}
we see that this a projective $\mc H (\mc R ,q )_\chi$-module.
Hence $P_n (V)$ also has this property.
It follows from \eqref{eq:2.7} that
\begin{equation}
P_n (V)^\Omega = R_\Omega \cdot P_n (V)
\end{equation}
is a direct summand of $P_n (V)$. \\
c) By a) and b) these cohomological dimensions are at most $r$.
On the other hand, we can easily find modules which do not have
projective resolutions of length smaller than $r$. Note that
\[
\mc A_\chi := \mc A \otimes_{Z(W)} \mh C_\chi \cong \mc O (T_\chi )
\]
where $T_\chi$ is the $r$-dimensional subtorus of $T$ consisting of
elements $t$ such that $t \big|_{Z(W)} = \chi$. Pick $t \in T_\chi$
and consider the parabolically induced module
\begin{equation}\label{eq:2.22}
I_t = \mr{Ind}_{\mc A}^{\mc H}(\mh C_t ) =
\mr{Ind}_{\mc A_\chi}^{\mc H (\mc R ,q )_\chi}(\mh C_t )
\end{equation}
With Theorem \ref{thm:2.12} we find
\begin{equation}\label{eq:2.23}
\mr{Ext}_{\mc H (\mc R ,q )_\chi}^r ( I_t ,I_t ) \cong
\mr{Ext}_{\mc A_\chi}^r (\mh C_t , I_t ) \cong
\mr{Ext}_{\mc O (T_\chi)}^r \Big( \mh C_t , \bigoplus_{w \in W_0}
\mh C_{w t} \Big) \cong \bigoplus_{w \in W_0 : w t = t} \mh C_{w t}
\end{equation}
Since this space is not 0, any resolution of $I_t$ by projective
$\mc H (\mc R ,q )_\chi$-modules has length at least $r$. This
calculation goes through in the bornological setting if we endow
all spaces with the fine bornology. $\qquad \Box$
\\[3mm]

For purposes of homological algebra it would be useful if we could
also construct projective resolutions for $\mc H$-modules that do not
admit a $Z(W)$-character. Unfortunately the authors do not know how
to achieve this in general. But we offer an alternative that comes
quite close. Let
\[
\mc H (\widetilde{\mc R} ,\tilde q) = \tilde G \ltimes \mc H (\mc R ,q)
\]
be a semisimple affine Hecke algebra as in \eqref{eq:2.20}.
Obviously $\mc H (\widetilde{\mc R} ,\tilde q)$ is a free
(left or right) $\mc H (\mc R ,q)$-module with basis $\{ T_g :
g \in \widetilde G \}$. Moreover for $(\pi ,V) \in \mr{Mod}(\mc H )$
the $\mc H (\widetilde{\mc R} ,\tilde q)$-module
\begin{equation}
\mr{Ind}_{\mc H (\mc R ,q)}^{\mc H (\widetilde{\mc R}
,\tilde q)} V = \mc H (\widetilde{\mc R} ,\tilde q) \otimes_{\mc H} V
\end{equation}
is isomorphic as an $\mc H$-module to
$\bigoplus_{g \in \widetilde G} V_g$, where the $\mc H$-module
structure on $V_g = (\pi_g ,V)$ is given by
\begin{equation}\label{eq:2.28}
\pi_g (h) \, v = \pi (T_g^{-1} h T_g ) \, v
\end{equation}
Clearly $V_g = V$ as an $\mc H (\mc R_{F_0},q)$-module. If $V$ admits
a $Z(W)$-character $\chi$ then $V_g$ differs only from $V$
in the sense that its $Z(W)$-character is $g \chi$.

Applying the construction of Corollary \ref{cor:2.11}.a) to
$\mc H (\widetilde{\mc R} ,\tilde q) \otimes_{\mc H} V $ as a
$\mc H (\widetilde{\mc R} ,\tilde q)$-module we get a resolution by
modules that are projective in Mod$ \big( \mc H (\widetilde{\mc R}
,\tilde q) \big)$ and in Mod$ \big( \mc H (\mc R ,q) \big)$.
In several cases this might be used to
find a resolution of $(\pi ,V)$ by projective $\mc H$-modules.

\begin{prop}\label{prop:2.dimh}
The cohomological dimensions of $\mr{Mod}(\mc H)$ and
$\mr{Mod}_{bor}(\mc H)$ are both equal to the rank of $X$.
\end{prop}
\emph{Proof.}
The cohomological dimension of Mod$ (\mc H )$ is the least number \\
$d \in \{0,1,2, \cdots, \infty \}$ such that
\[
\mr{Ext}^n_{\mc H} (U,V) = 0 \qquad
\forall \, U,V \in \mr{Mod}(\mc H) \,, \forall n > d
\]
Let $t \in T$ and consider the module $I_t = \mr{Ind}_{\mc A}^{\mc H}
(\mh C_t )$. In view of Theorem \ref{thm:2.12}
\[
\mr{Ext}_{\mc H}^{\mr{rk}(X)} ( I_t ,I_t ) \cong
\mr{Ext}_{\mc A}^{\mr{rk}(X)} (\mh C_t , I_t ) \cong
\mr{Ext}_{\mc O (T)}^{\mr{rk}(X)} \Big( \mh C_t , \bigoplus_{w \in W_0}
\mh C_{w t} \Big) \cong \bigoplus_{w \in W_0 : w t = t} \mh C_{w t}
\]
Therefore $d \geq \mr{rk}(X)$. This argument also works in
$\mr{Mod}_{bor} (\mc H )$, provided that we endow all spaces with the
fine bornology.

On the other hand, let $U,V \in \mr{Mod} (\mc H )$ be arbitrary and
consider the $\mc H (\widetilde{\mc R} ,\tilde q)$-modules
$\mr{Ind}_{\mc H}^{\mc H (\widetilde{\mc R} ,\tilde q)}(U)$ and
$\mr{Ind}_{\mc H}^{\mc H (\widetilde{\mc R} ,\tilde q)}(V)$.
\begin{equation}\label{eq:2.26}
\begin{array}{lll}
\mr{Ext}_{\mc H}^n (U,V) & \subset &  \bigoplus_{g \in \tilde G}
\mr{Ext}_{\mc H}^n \big( U,V_g \big) \: \cong \: \mr{Ext}_{\mc H}^n
\big( U, \mr{Ind}_{\mc H}^{\mc H (\widetilde{\mc R} ,\tilde q)}(V)
\big) \cong \\
& & \mr{Ext}_{\mc H (\widetilde{\mc R} ,\tilde q)}^n \big(
\mr{Ind}_{\mc H}^{\mc H (\widetilde{\mc R} ,\tilde q)}(U)
,\mr{Ind}_{\mc H}^{\mc H (\widetilde{\mc R} ,\tilde q)}(V) \big)
\end{array}
\end{equation}
Assume $n > \mr{rk} (X)$. According to Corollary \ref{cor:2.11}.c)
the cohomological dimension of Mod$( \mc H (\widetilde{\mc R} ,\tilde q))$
is rk$ (X)$, so right hand side of \eqref{eq:2.26} is 0. Hence
$\mr{Ext}_{\mc H}^n (U,V) = 0$ and we conclude that $d \leq \mr{rk}(X)$.
The same reasoning shows that the cohomological dimension of
$\mr{Mod}_{bor}(\mc H)$ is rk$ (X). \qquad \Box$
\\[3mm]
Recall that a resolution $(P_* ,d_* )$ of a module $V$ is of finite
type if all the modules $P_n$ are finitely generated, and moreover
$P_n = 0$ for all $n$ larger then some number.

\begin{cor}\label{cor:2.ftp}
Let $V$ be a finitely generated $\mc H$-module. Then $V$ admits a
finite type projective resolution.
\end{cor}
\emph{Proof.} Because $\mc H$ is Noetherian, every submodule of a
finitely generated $\mc H$-module is itself finitely generated.

By assumption there exist a surjective $\mc H$-module map $d_0 : \mc
H^{m_0} \to V$, for some $m_0 \in \mh N$. Then ker $d_0$ is again
finitely generated, so we can find a surjection $d_1 : \mc H^{m_1}
\to \ker d_0$. Continuing this process we construct a resolution
$(P_n = \mc H^{m_n} ,d_n )$ of $V$, consisting of free $\mc
H$-modules of finite rank. Because the global dimension of $\mc H$
is rk$(X)$, the module ker $d_n$ must be projective $\forall n \geq
\mr{rk}(X) - 1$ \cite[Proposition VI.2.1]{CaEi}. Hence
\[
0 \leftarrow V \xleftarrow{d_0} P_0 \xleftarrow{d_1} \cdots
\xleftarrow{d_{n-1}} P_{n-1} \leftarrow \ker d_{n-1} \leftarrow 0
\]
is a finite type projective resolution of $V. \qquad \Box$
\\[4mm]

\section{Projective resolutions for Schwartz algebras}
\label{sec:2.3}

We will show that all the resolutions from the previous section
can be induced from $\mc H$ to $\mc S$. Most importantly, we will
construct a projective bimodule resolution of $\mc S$. This
requires that we complete the $\mc H$-modules to Fr\'echet
$\mc S$-modules.
A convenient technique to achieve this in great generality
is with completed bornological tensor products, and this is
the viewpoint we chose to take in this section. However,
for finite dimensional tempered modules it is not necessary
to use bornologies. See the remark after Corollary \ref{cor:2.9}.

Let $V \in \mr{Mod}_{bor}(\mc S )$. According to
\cite[Theorem 42]{Mey1} we have
\begin{equation}\label{eq:2.szw}
\mc S (Z(W)) \hot_{\mh C [Z(W)]} V =
\mc S (Z(W)) \hot_{\mc S (Z(W))} V
\end{equation}
If $V$ has finite dimension, then \eqref{eq:2.szw} also holds with
algebraic tensor products. The reader is invited to check this, by
reduction to the case where $V$ admits a unique $Z(W)$-character.

Because $\mc H$ is a free $\mc H (W_f ,q) \otimes \mh C [Z(W)]$-module,
both algebraically and with the fine bornology, we have
\begin{equation}
\mc H \hot_{\mc H (W_f ,q) \otimes \mh C [Z(W)]} V =
\mc H \otimes_{\mc H (W_f ,q) \otimes \mh C [Z(W)]} V
\end{equation}
So if we induce $P_n (V )$ from $\mc H$ to $\mc S$ in
the bornological fashion we get the module
\begin{equation}\label{eq:2.pntv}
\begin{split}
P_n^t (V) \: := \: \mc S \hot_{\mc H} P_n (V) \: = \: \mc S
\hot_{\mc H} \bigoplus_{f : \dim f = n} \mc H \hot_{\mc H (W_f ,q)
\otimes \mh C [Z(W)]} V \otimes_{\mh C} \mh C \{ f \} \\
= \: \bigoplus_{f : \dim f = n} \mc S \hot_{\mc H (W_f ,q) \otimes
\mc S (Z(W))} V \otimes_{\mh C} \mh C \{ f \}
\end{split}
\end{equation}
The maps $d_n : P_n (V) \to P_{n-1}(V)$ extend naturally to
\[
d_n^t : P_n^t (V) \to P_{n-1}^t (V)
\]
The action of $\Omega$ on $P_n (V)$ also extends to $P_n^t (V)$,
so we can construct $P_n^t (V )^\Omega$. By \eqref{eq:2.7}
\begin{equation}\label{eq:2.ro}
P_n^t (V )^\Omega = R_\Omega \cdot P_n^t (V)
\end{equation}
is a direct summand of $P_n^t (V)$. Clearly $P_n^t (V)$
and $P_n^t (V )^\Omega$ are finitely generated $\mc S$-modules
if $V$ has finite dimension.

We consider the important case $V = \mc S$. The topology and the
bornology on $\mc S$ give rise to a topology and a bornology on
$P_n^t (\mc S )$. For $n,m,k \in \mh N \,,
f \subset \overline{A_\es}$ we have the continuous seminorms
\begin{align*}
& p_{m,k,f} : \mc S \hot_{\mc H (W_f ,q) \otimes \mc S (Z(W))}
\mc S \otimes_{\mh C} \mh C \{ f \} \to [0,\infty ) \\
& p_{m,k,f}(y) = \inf \Big\{ \sum_i p_m (h_i ) p_k (h'_i) :
\sum_i h_i \otimes h'_i \otimes f = y \Big\}
\end{align*}
which define a Fr\'echet topology on this space. The topology on
$P_n^t (\mc S )$ is defined by the norms $p_{m,k} := \sum_f p_{m,k,f}$.
We endow these modules with the precompact bornology.
We note that $d_n^t$ is continuous and bounded and that $P_n (\mc S )$
is dense in $P_n^t (\mc S )$.

In view of \eqref{eq:2.ef} we have
\[
P_n^t (\mc S )^\Omega = \bigoplus_{f : dim f = n} \mc S (\mc R_{F_0},q)
\hot_{\mc H (W_f ,q)} \mc S (\mc R ,q) \cong \bigoplus_{f : dim f = n}
\Big( \mc S (\mc R_{F_0},q) \hot_{\mh C} \mc S (\mc R ,q)^{op} \Big) e_f
\]
Using Lemma \ref{lem:2.10} and Theorem \ref{thm:2.7} both for
$\mc S (\mc R_{F_0},q)$ and for $\mc S (\mc R ,q)$ we see that
there is a number $C_{m,k,f} > 0$ such that
\begin{equation}\label{eq:2.est}
\begin{array}{ll}
\sum_{w \in W^\af ,w' \in W} |h_{w,w'}| (\mc N (w) +1 )^m
(\mc N (w') +1 )^k & \leq \\
C_b^2 \sup_{w \in W^\af ,w' \in W} |h_{w,w'}|
(\mc N (w) +1 )^{m+b} (\mc N (w') +1 )^{k+b} & \leq \\
C_b^2 p_{m+b,k+b} \Big( \sum_{w \in W^\af ,w' \in W} h_{w,w'}
(N_w \otimes N_w' ) e_f  \otimes f \Big) & \leq \\
C_{m,k,f} p_{m+2b,k+2b,f} \Big( \sum_{w \in W^\af} \sum_{w' \in W}
h_{w,w'} N_w \otimes N_{w'} \Big) &
\end{array}
\end{equation}

\begin{thm}\label{thm:2.8}
Consider $\mc S$ as a $\mc S$-bimodule.
\begin{equation}\label{eq:2.21}
0 \longleftarrow \mc S \xleftarrow{\; d_0^t \;}
P_0^t (\mc S )^\Omega \xleftarrow{\; d_1^t \;}
P_1^t (\mc S )^\Omega \longleftarrow \cdots
\xleftarrow{\; d_r^t \;} P_r^t (\mc S )^\Omega \longleftarrow 0
\end{equation}
is a $\mc S \hot \mc S^{op}$-module resolution of $\mc S$, with a
continuous bounded contraction. Every $P_n^t (\mc S )$ is a bornologically
projective $\mc S$-module, both from the left and from the right.
Moreover if $\mc R$ is semisimple then $P_n^t (\mc S )$ is also
projective as a $\mc S \hot \mc S^{op}$-module.
\end{thm}
\emph{Proof.}
To show that the differential complex $\big( P_*^t (\mc S)^\Omega ,
d_*^t \big)$ is contractible we use Proposition \ref{prop:2.3} and
Theorem \ref{thm:2.6}. The composition of \eqref{eq:2.8} with
\eqref{eq:2.9} is the bijection
\begin{align*}
& \tilde \phi : C_* (\Sigma ) \otimes_{\mh C} \mc S \to P_* (\mc S)^\Omega \\
& \tilde \phi (\sigma \otimes h' ) = R_\Omega \big( T_w
\otimes_{\mc H (W_f ,q) \otimes \mh C [Z(W)]} T_w^{-1} h' \otimes f \big)
\end{align*}
where $\sigma = w \overline f$ with $w \in W^\af$. Let $\gamma$ be as
in Proposition \ref{prop:2.3}. We claim that
\[
\tilde \gamma := \tilde \phi (\gamma \otimes \mr{id}_\mc S )
{\tilde \phi}^{-1}
\]
extends continuously to the required contraction.
Suppose that $w' \in W ,\, w \in W^\af \cap w' \Omega$ and
$\sigma = w' \overline f = w \overline f$. Then we have explicitly
\begin{equation}\label{eq:2.10}
\tilde \phi \big( R_\Omega (N_{w'} \otimes_{\mc H (W_f ,q)
\otimes \mh C [Z(W)]} h' \otimes f ) \big) =
\tilde \phi \big( \gamma (\sigma ) \otimes N_w h' \big) =
\tilde \phi \Big( \sum_\tau \gamma_{\sigma \tau} \tau \otimes N_w h' \Big)
\end{equation}
By Lemma \ref{lem:2.2} and condition 3) of Proposition \ref{prop:2.3}
the coefficient $\gamma_{\sigma \tau}$ can only be nonzero if there
exist $u \leq_A w$ and a facet $f'$ of $A_\es$ such that
$\tau = u \overline{f'}$. This crucial for the following estimates.
For every relevant $\tau$ we pick such a $u \in W^\af$ and we write
(a little sloppily) $\gamma_{w u} = \gamma_{\sigma \tau}$. Then
\eqref{eq:2.10} equals
\begin{equation}
\begin{array}{ll}
\tilde \phi \Big( \sum_{f'} \sum_{u \in W^\af : u \leq_A w}
\gamma_{w u} (u \overline{f'}) \otimes N_w h' \Big) & = \\
R_\Omega \Big( \sum_{f'} \sum_{u \in W^\af : u \leq_A w}
\gamma_{w u} N_u \otimes_{\mc H (W_{f'} ,q) \otimes \mh C [Z(W)]}
N_u^{-1} N_w h' \otimes f' \Big) & = \\
\sum_{f'} \sum_{u \in W^\af : u \leq_A w} R_\Omega \big(
\gamma_{wu} N_u \otimes_{\mc H (W_{f'} ,q) \otimes \mh C [Z(W)]}
N_{u^{-1} w} h' \otimes f' \big) &
\end{array}
\end{equation}
Notice that we used $u \leq_A w$ in the last step. Every element of
$P_n^t (\mc S)^\Omega$ can be written as a finite sum (over facets $f$)
of elements of the form
\[
R_\Omega \, y = R_\Omega \sum_{w \in W^\af} \sum_{w' \in W} h_{w,w'}
N_w \otimes_{\mc H (W_{f'} ,q) \otimes \mc S (Z(W) )} N_{w'} \otimes f
\]
with $(h_{w,w'}) \in \mc S (W^\af \times W)$.
According to the above calculation
\[
\tilde \gamma (R_\Omega \, y) = R_\Omega \sum_{f'} \sum_{w' \in W}
\sum_{u,w \in W^\af : u \leq_A w} \gamma_{wu} h_{w,w'} N_u
\otimes_{\mc H (W_{f'} ,q) \otimes \mc S (Z(W))} N_{u^{-1} w} N_{w'}
\otimes f'
\]
Using (in this order) condition 4) of Proposition \ref{prop:2.3}, Theorem
\ref{thm:2.7}, Lemma \ref{lem:2.1} and \eqref{eq:2.est} we estimate
\begin{align*}
& p_{m,k} \Big( \sum_{w' \in W} \sum_{u,w \in W^\af : u \leq_A w}
\gamma_{wu} h_w N_u \otimes_{\mc H (W_{f'} ,q) \otimes \mc S (Z(W) )}
N_{u^{-1} w} h' \otimes f' \Big) \quad & \leq \\
& \sum_{w' \in W} \sum_{w \in W^\af} M_\gamma |h_{w,w'} | p_m
\big( \sum_{u \in W^\af : u \leq_A w} N_u \big) p_k \big( N_{u^{-1} w}
N_{w'} \big) \quad & \leq \\
& M_\gamma \sum_{w' \in W} \sum_{w \in W^\af} |h_{w,w'} | (\mc N (w) + 1)^m
C_q (\mc N (w) + 1)^{k + b} (\mc N (w') +1 )^{k+b} \quad & \leq \\
& M_\gamma C_q C_{k+m+2b,k+b,f} p_{k+m+3b , k+2b} (y)
\end{align*}
Since $R_\Omega$ is a continuous operator on $P_n^t (\mc S)$
it follows that $\tilde \gamma$ is well-defined and continuous on
$P_n^t (\mc S)^\Omega$. Since $P_n^t (\mc S )$ carries the precompact
bornology $\tilde \gamma$ is automatically bounded. Moreover
\[
\tilde \phi (\delta_n \otimes \mr{id}_\mc S ) {\tilde \phi}^{-1} = d_n
\]
so condition 1) of Proposition \ref{prop:2.3} assures that
\begin{equation}\label{eq:2.17}
\tilde \gamma d^t + d^t \tilde \gamma = \text{ id}
\end{equation}
on $P_* (\mc S)^\Omega$. Because $P_* (\mc S)^\Omega$ is dense in
$P_*^t (\mc S)^\Omega$ and the maps in \eqref{eq:2.17} are continuous,
this relation holds on the whole of $P_*^t (\mc S)^\Omega$. So the
differential complex $\big( P^t_* (\mc S)^\Omega , d^t_* \big)$
indeed has a bounded contraction.

For any facet $f$ the space $\mc S$ is a bornologically free
$\mc H (W_f ,q) \otimes \mc S (Z(W))$-module. Hence $P_n^t (\mc S)$
is a bornologically projective $\mc S$-module, both from the left and
from the right. If $\mc R$ is semisimple then by \eqref{eq:2.ef}
$P_n^t (V)$ is direct sum of bimodules of the form
$(\mc S \hot \mc S^{op}) e_f$. Hence $P_n^t (V)$ is
$\mc S \hot \mc S^{op}$-projective.

By \eqref{eq:2.ro} $P_n^t (\mc S )^\Omega$ enjoys the same
projectivity properties. $\qquad \Box$
\\[2mm]

\begin{cor}\label{cor:2.9}
\begin{description}
\item[a)] Let $V$ be any bornological $\mc S$-module.
\[
0 \longleftarrow V \xleftarrow{\; d_0^t \;} P_0^t (V )^\Omega
\xleftarrow{\; d_1^t \;} P_1^t (V )^\Omega \longleftarrow \cdots
\xleftarrow{\; d_r^t \;} P_r^t (V )^\Omega \longleftarrow 0
\]
is a bornological resolution of $V$.
\item[b)] If $V$ admits the $Z(W)$-character $\chi$ then every
module $P_n^t (V)^\Omega$ is projective in
$\mr{Mod}_{bor} \big(\mc S (\mc R ,q)_\chi \big)$.
\item[c)] If moreover $V$ has finite dimension then $P_n^t (V)^\Omega$
is also projective in $\mr{Mod} \big(\mc S (\mc R ,q)_\chi \big)$.
\end{description}
\end{cor}
\emph{Proof.} a)
Apply $\otimes_{\mc S} V$ to \eqref{eq:2.21} and use the projectivity
of $P_n^t (\mc S )^\Omega$ as a right $\mc S$-module.\\
b) From Corollary \ref{cor:2.11}.b) we know that $P_n (V )^\Omega$ is
projective in $Mod_{bor} \big( \mc H (\mc R ,q )_\chi \big)$ so
\[
P_n^t (V )^\Omega \cong \mc S (\mc R ,q )_\chi
\hot_{\mc H (\mc R ,q )_\chi} P_n (V )^\Omega
\]
is projective in $\mr{Mod}_{bor} \big( \mc S (\mc R ,q )_\chi \big)$. \\
c) For any facet $f$
\[
\mc S \hot_{\mc H (W_f,q) \otimes \mc S (Z(W))} V =
\mc S (\mc R ,q )_\chi \hot_{\mc H (W_f ,q)} V =
\mc S (\mc R ,q )_\chi \otimes_{\mc H (W_f ,q)} V =
\mr{Ind}_{\mc H (W_f ,q)}^{\mc S (\mc R ,q )_\chi} V
\]
is a projective $\mc S (\mc R ,q )_\chi$-module. In view
of \eqref{eq:2.ro} this implies that $P_n^t (V)$ and $P_n^t (V )^\Omega$
are also projective in Mod$ \big( \mc S (\mc R ,q )_\chi \big).
\qquad \Box$
\\[3mm]
\textbf{Remark.}\\
If $V$ is a finite dimensional tempered module with $Z(W)$-character
$\chi$ then the proof of Corollary \ref{cor:2.9} does not rely on the
properties of bornology. Indeed, in this situation we may simply use the
algebraic tensor product in the definition of $P^t_n (V)$, since the
algebraic tensor product is already complete as a locally convex
vector space. The continuity proof of the contraction is analogous
to and in fact somewhat simpler than the above proof for the case
$V = \mc S$. Hence the algebraic tensor product of the resolution of
Corollary 2.3 a) by $\mc S(R,q)_\chi$ yields the resolution of Corollary
\ref{cor:2.9}.a).
\vspace{4mm}

\section{Isocohomological inclusions}

We will show that the inclusion $\mc H \to \mc S$ is isocohomological.
As an intermediate step we do the same for algebras and modules
corresponding to a fixed $Z(W)$-character.

Similar results for Schwartz algebras of reductive $p$-adic groups
were proven by Meyer \cite[Theorems 21, 27 and 29]{Mey2} with highly
sophisticated techniques. Maybe our bounded contraction from Section
\ref{sec:2.1} can be used to simplify these proofs.

\begin{thm}\label{thm:dimsx}
Let $\chi$ be a unitary $Z(W)$-character.
\item[a)] The inclusion $\mc H (\mc R ,q)_\chi \to
\mc S (\mc R ,q )_\chi$ is isocohomological.
\item[b)] The cohomological dimension of $\mr{Mod}_{bor} \big(
\mc S (\mc R ,q )_\chi \big)$ equals $r = \mr{rk} (R_0 )$.
\end{thm}
\emph{Proof.} a)
From \eqref{eq:2.ind} and \eqref{eq:2.pntv} it follows that
\[
\begin{array}{lll}
P_n (\mc H (\mc R ,q )_\chi ) & \cong & \bigoplus\limits_{f : \dim f = n}
\mc H (\mc R ,q )_\chi \hot_{\mc H (W_f ,q)}
\mc H (\mc R ,q )_\chi \otimes_{\mh C} \mh C \{ f \} \\
P_n^t (\mc S (\mc R ,q )_\chi ) & \cong & \bigoplus\limits_{f : \dim f = n}
\mc S (\mc R ,q )_\chi \hot_{\mc H (W_f ,q)}
\mc S (\mc R ,q )_\chi \otimes_{\mh C} \mh C \{ f \}
\end{array}
\]
Exactly as in the proof of Theorem \ref{thm:2.6} we can see that these
are projective as bornological bimodules for $\mc H (\mc R ,q )_\chi$
respectively $\mc S (\mc R ,q )_\chi$. In view of \eqref{eq:2.7} and
\eqref{eq:2.ro} the same holds for $P_n (\mc H (\mc R ,q )_\chi )^\Omega$
and $P_n^t (\mc S (\mc R ,q )_\chi )^\Omega$. Combined with Corollaries
\ref{cor:2.11}.a) and \ref{cor:2.9}.a) this yields condition 1) of
Theorem \ref{thm:2.15}.\\
b) By Corollary \ref{cor:2.9} the cohomological dimension of
$\mr{Mod}_{bor} \big( \mc S (\mc R ,q )_\chi \big)$ is at most
$r = \mr{rk} (R_0 )$. If $t \in T$ is unitary then by
\cite[Proposition 4.19]{Opd1} the module $I_t$ from \eqref{eq:2.22}
is tempered. Together with \eqref{eq:2.23} this gives
\[
\mr{Ext}_{\mc S (\mc R ,q )_\chi}^r (I_t ,I_t ) \cong
\mr{Ext}_{\mc H (\mc R ,q )_\chi}^r (I_t ,I_t ) \neq 0
\]
Hence this cohomological dimension is at least $r. \qquad \Box$
\\[3mm]

\begin{thm}\label{thm:2.16}
\begin{description}
\item[a)] The inclusion $\mc H \to \mc S$ is isocohomological.
\item[b)] The cohomological dimension of $\mr{Mod}_{bor}(\mc S )$
equals the rank of $X$.
\end{description}
\end{thm}
\emph{Proof.} a)
Let $(\widetilde{\mc R} ,\tilde q)$ be as in \eqref{eq:2.20}.
Recall that
\[
\begin{array}{lll}
\mc H (\widetilde{\mc R} ,\tilde q) & \cong &
\tilde G \ltimes \mc H (\mc R ,q) \\
\mc S (\widetilde{\mc R} ,\tilde q) & \cong &
\tilde G \ltimes \mc S (\mc R ,q)
\end{array}
\]
We know from Theorem \ref{thm:dimsx}.a) that the inclusion
$\mc H (\widetilde{\mc R} ,\tilde q) \to \mc S (\widetilde{\mc R}
,\tilde q)$ is isocohomological. Therefore we can use an argument
from the proof of \cite[Theorem 58]{Mey1}. The functor
\begin{equation}\label{eq:2.13}
\mr{Mod} (B ) \to \mr{Mod} (\tilde G \ltimes B ) \,:\; V \to
\mr{Ind}_B^{\tilde G \ltimes B} (V) = (\tilde G \ltimes B ) \otimes_B V
\end{equation}
is exact for any $\tilde G$-algebra $B$. Hence in $\mr{Der}_{bor}
(\tilde G \ltimes \mc S )$ we have
\begin{equation}\label{eq:2.15}
\begin{array}{lll}
\tilde G \ltimes \mc S & \cong & (\tilde G \ltimes \mc S )
\hot_{\tilde G \ltimes \mc S}^{\mh L} (\tilde G \ltimes \mc S ) \\
& \cong & (\tilde G \ltimes \mc S )
\hot_{\tilde G \ltimes \mc H}^{\mh L} (\tilde G \ltimes \mc S ) \\
& \cong & (\tilde G \ltimes \mc S ) \hot_{\tilde G \ltimes \mc H
}^{\mh L} \tilde G \ltimes \mc H \hot_{\mc H}^{\mh L} \mc S \\
& \cong & (\tilde G \ltimes \mc S ) \hot_{\mc H}^{\mh L} \mc S \\
& \cong & \mr{Ind}_{\mc S}^{\tilde G \ltimes \mc S}
\big( \mc S \hot_{\mc H}^{\mh L} \mc S \big)
\end{array}
\end{equation}
We want to show that this implies condition 2) of Theorem
\ref{thm:2.15} for the inclusion $\mc H \to \mc S$. However we
have to be a little careful, as the functor \eqref{eq:2.13} is not
injective on objects. Namely, $\mc H$-modules like $V$ and $V_g$
in \eqref{eq:2.28}, which are conjugate by an element of $\tilde
G$, have the same image under \eqref{eq:2.13}. It follows from
\eqref{eq:2.15} that
\begin{equation}\label{eq:2.gss}
\mh C [\tilde G ] \otimes_{\mh C} \mr{Tor}^{\mc H}_n (\mc S ,\mc S )
\cong \left\{ \begin{array}{lll}
\tilde G \ltimes \mc S & \mr{if} & n = 0 \\
0 & \mr{if} & n > 0
\end{array} \right.
\end{equation}
Obviously the multiplication map
\[
\mr{Tor}^{\mc H}_0 (\mc S ,\mc S ) \cong \mc S \hot_{\mc H} \mc S \to \mc S
\]
is surjective. In view of \eqref{eq:2.gss} it must
also be injective, and therefore
\[
\mr{Tor}^{\mc H}_n (\mc S ,\mc S )
\cong \left\{ \begin{array}{lll}
\mc S & \mr{if} & n = 0 \\
0 & \mr{if} & n > 0
\end{array} \right.
\]
Let
\begin{equation}\label{eq:2.spp}
0 \leftarrow \mc S \leftarrow P_0 \leftarrow P_1
\leftarrow \cdots
\end{equation}
be a bornological resolution of $\mc S$ by projective $\mc H$-modules.
We already know that the homology of \eqref{eq:2.spp} vanishes
in all degrees.
Moreover $\mr{Ind}_{\mc H}^{\tilde G \ltimes \mc H} ( P_* )$ is a
resolution of $\tilde G \ltimes \mc H$. Theorems \ref{thm:dimsx}.a)
and \ref{thm:2.15} assure that the differential
complex $\mr{Ind}_{\mc H}^{\tilde G \ltimes \mc H} \big( \mc S
\hot_{\mc H} P_* \big)$ is a bornological resolution of $\tilde G
\ltimes \mc S$. In particular it admits a bounded $\mh C$-linear
contraction. Hence $\mc S \hot_{\mc H} P_* $ also admits a bounded
contraction, i.e. it is an exact sequence in $\mr{Mod}_{bor} (\mc S )$.
This shows that the natural map
\begin{equation}
\mc S \hot_{\mc H}^{\mh L} \mc S \to
\mc S \hot_{\mc S}^{\mh L} \mc S
\end{equation}
is an isomorphism. We conclude that $\mc H \to \mc S$ is indeed
isocohomological.\\
b)
In view of part a) and Proposition \ref{prop:2.dimh} the
cohomological dimension of $\mr{Mod}_{bor} (\mc S )$ is at
most $\mr{rk}(X)$. If $t \in T$ is unitary then by
\cite[Proposition 4.19]{Opd1} the module $I_t$ from \eqref{eq:2.22}
is tempered. From a) and the proof of Proposition \ref{prop:2.dimh}
we see that
\[
\mr{Ext}_{\mc S}^{\mr{rk}(X)} (I_t ,I_t ) \cong
\mr{Ext}_{\mc H}^{\mr{rk}(X)} (I_t ,I_t ) \neq 0
\]
Hence this cohomological dimension is at least $\mr{rk}(X).
\qquad \Box$
\\[3mm]
\textbf{Remark.}\\
In the same way one can show that the cohomological dimension of
the category $\mr{Mod}_{Fr\acute e} (\mc S )$ of continuous Fr\'echet
$\mc S$-modules is the rank of $X$. To make this a
meaningful statement we make this into an exact category as follows.

All morphisms are required to be continuous and $\hot$ is the
completed projective tensor product. Only extensions and resolutions
that admit a continuous $\mh C$-linear splitting are called exact.
This category has enough projective objects and has countable
projective limits. However it does neither have enough injective
objects nor inductive limits.

%% file: tempered3.tex
\chapter{The Euler-Poincar\'e characteristic}

\section{Elliptic representation theory}

Elliptic representation theory is a general notion that can be
developed for many groups and algebras
\cite{Art,Kaz,Ree,ScSt,Wal}. The idea is that one considers all
virtual representations of an algebra, modulo those that are
induced from certain specified subalgebras. This should yield
interesting equivalence classes of representations if the
subalgebras are chosen cleverly.

For example in a reductive $p$-adic group one can consider
the collection of proper parabolic subgroups. The resulting
space of representations contains among others all square
integrable representation. It can be studied by means of
certain integrals over the regular elliptic conjugacy classes,
cf. \cite{Kaz,Bez,ScSt}.

In the context of the elliptic representation theory for
Iwahori-spherical representations of a p-adic Chevalley group
Reeder [Ree] was led to the general definition of elliptic
representation theory for a finite group relative to a given
representation. Let $(\rho ,E)$ be a real representation of a
finite group $\Gamma$. We define an elliptic pairing on
$\mr{Mod}_{fin}(\Gamma )$ by
\begin{equation}
e_\Gamma (U,V) := \sum_{n=0}^\infty (-1 )^n \dim
\mr{Hom}_\Gamma \big( U \otimes \wig^n E ,V \big)
\end{equation}
We call an element $\gamma \in \Gamma$ elliptic (with respect
to $E$) if $E^{\rho (\gamma)} = 0$. Since this property is preserved
under conjugation, we can use the same terminology for
conjugacy classes. Let $\mc L$ be the set of subgroups
$H \subset \Gamma$ such that $E^{\rho (H)} \neq 0$. The space of
elliptic trace functions on $\Gamma$ is defined as
\begin{equation}
Ell (\Gamma ) := G_{\mh C}(\Gamma ) \Big/ \sum_{H \in \mc L}
\mr{Ind}_H^\Gamma \big( G_{\mh C}(H) \big)
\end{equation}

\begin{thm}\label{thm:3.8}
\begin{description}
\item[a)] The dimension of $Ell (\Gamma )$ equals the number of
elliptic conjugacy classes of $\Gamma$.
\item[b)] $e_\Gamma$ induces a Hermitian inner product
on $Ell (\Gamma )$.
\item[c)] For all $\chi, \chi' \in G_{\mh C} (\Gamma )$ we have
\[
e_\Gamma (\chi,\chi' ) = \sum_{\gamma \in \Gamma}
\frac{\det \, (\mr{id}_E - \rho (\gamma ) )}{|\Gamma|}
\overline{\chi (\gamma)} \chi' (\gamma)
\]
\end{description}
\end{thm}
\emph{Proof.} See \cite[\S 2]{Ree}.
$\qquad \Box$ \\[3mm]

Assume now that $X$ is a lattice in $E$ (so
$E = X \otimes_{\mh Z} \mh R$) which is stable under the action
of $\Gamma$. We will show that Theorem \ref{thm:3.8} can be
generalized to the group $\Gamma \ltimes X$. Of course affine
Weyl groups are important examples of such groups.

In what follows expression like $\gamma x$ always should be
interpreted as the product in $\Gamma \ltimes X$. If we want
to make $\gamma$ act on $x$ then we write $\rho (\gamma ) x$.
We extend this to an action of $\Gamma \ltimes X$ on $X$ by
\[
\rho (y \gamma ) x = y + \rho (\gamma) x
\]
Let $t \in T = \mr{Hom}_{\mh Z} (X,\mh C^\times )$. It is known
from \cite{Cli} that there is a natural bijection between
irreducible representations of $\Gamma_t = \{ \gamma \in
\Gamma : t \circ \rho (\gamma ) = t \}$ and irreducible
representations of $\Gamma \ltimes X$ with central character
$\Gamma t \in T / \Gamma$.
It is given explicitly by
\begin{equation}\label{eq:3.2}
\mr{Ind}_t : V \to
\mr{Ind}_{\Gamma_t \ltimes X}^{\Gamma \ltimes X} V_t
\end{equation}
where $V_t$ means that we regard $V$ as a $X$-representation
with character $t$.

We call an element $\gamma x \in \Gamma \ltimes X$ elliptic
if it has an isolated fixpoint in $E$. It is easily seen that
this is the case if and only if $\gamma \in \Gamma$ is
elliptic. We have
\[
x \, y \gamma \, (-x) =
(x - \rho (\gamma) x ) \, y \gamma \in \Gamma \ltimes X
\]
so all elements of $(y + (\mr{id}_E - \rho (\gamma )) X) \gamma$
are conjugate in $\Gamma \ltimes X$. If $\gamma$ is elliptic then
the lattice $(\mr{id}_E - \rho (\gamma ) )X$ is of finite
index in $X$. Consequently there are only finitely many elliptic
conjugacy classes in $\Gamma \ltimes X$.

Let $U$ and $V$ be $\Gamma \ltimes X$ modules of finite length,
i.e. finite dimensional. We define the Euler-Poincar\'e
characteristic
\begin{equation}
EP_{\Gamma \ltimes X} (U,V) = \sum_{n=0}^\infty (-1 )^n \dim
\mr{Ext}^n_{\Gamma \ltimes X} (U,V)
\end{equation}
This kind of pairing stems from Schneider and Stuhler
\cite[\S III.4]{ScSt}, who studied it for reductive $p$-adic
groups. The space of elliptic trace functions on $\Gamma
\ltimes X$ is
\begin{equation}
Ell (\Gamma \ltimes X) := G_{\mh C}(\Gamma \ltimes X) \Big/
\sum_{H \in \mc L} \mr{Ind}_{H \ltimes X}^{\Gamma \ltimes X}
\big( G_{\mh C}(H \ltimes X) \big)
\end{equation}
For every $t \in T$ we consider the elliptic representation
theory of $\Gamma_t$ with respect to the cotangent space
to $T$ at $t$. We note that $\mr{Ind}_t$ induces a map
$Ell (\Gamma_t ) \to Ell (\Gamma \ltimes X)$. Let $H_{ell}$
denote the set of elliptic elements in a group $H$, and let
$\sim_H$ be the equivalence relation "conjugate by an
element of $H$".

\begin{thm}\label{thm:3.9}
\begin{description}
\item[a)] The dimension of $Ell (\Gamma \ltimes X)$ equals the
number of elliptic conjugacy classes of $\Gamma \ltimes X$.
\item[b)] $EP_{\Gamma \ltimes X}$ induces a Hermitian inner
product on $Ell (\Gamma \ltimes X)$.
\item[c)] The map $\mr{Ind}_t : Ell (\Gamma_t) \to Ell
(\Gamma \ltimes X)$ induced by \eqref{eq:3.2} is an isometry:
\[
EP_{\Gamma \ltimes X} (\mr{Ind}_t U , \mr{Ind}_t V ) =
e_{\Gamma_t} (U,V)
\]
for all finite dimensional $\Gamma_t$-representations and
$U$ and $V$.
\item[d)] The map $\bigoplus_{t \in T / \Gamma} \mr{Ind}_t
: \bigoplus_{t \in T / \Gamma} Ell (\Gamma_t)
\to Ell (\Gamma \ltimes X)$ is an isomorphism.
\end{description}
\end{thm}
\emph{Proof.}
For $U,V$ and $t$ as above we have by Frobenius reciprocity
\begin{equation}\label{eq:3.4}
\mr{Ext}^n_{\Gamma \ltimes X} (\mr{Ind}_t U , \mr{Ind}_t V )
\cong \mr{Ext}^n_{\Gamma_t \ltimes X} (U_t ,
\mr{Ind}_{\Gamma_t \ltimes X}^{\Gamma \ltimes X} V_t )
\end{equation}
Because two $\Gamma_t \ltimes X$-representations with different
central characters admit only trivial extensions, \eqref{eq:3.4}
is isomorphic to $\mr{Ext}_{\Gamma_t \ltimes X}^n (U_t ,V_t )$.
Inside the group algebra
\[
\mc A := \mh C [X] \cong \mc O (T)
\]
we have the ideal of functions vanishing at $t \in T$:
\[
I_t := \{ f \in \mc A : f(t) = 0 \}
\]
Let us denote the completion of $\mc A$ with respect to the
powers of this ideal by $\hat{\mc A_t}$. Clearly
\[
(\Gamma_t \ltimes \hat{\mc A_t} )
\otimes_{\Gamma_t \ltimes X} U_t = U_t
\]
as $\Gamma_t \ltimes X$-modules. Completing is an exact functor
so \eqref{eq:3.4} becomes
\begin{equation}
\mr{Ext}^n_{\mh C [\Gamma_t \ltimes X]} (U_t , V_t ) \cong
\mr{Ext}^n_{\Gamma_t \ltimes \hat{\mc A_t}} (U_t , V_t )
\end{equation}
Because the $\Gamma_t$-module $I_t^2$ has finite codimension in
$\mc A$ there exists a $\Gamma_t$-module $E_t \subset \mc A$ such that
\begin{equation}
\mc A = \mh C \oplus E_t \oplus I_t^2
\end{equation}
As a $\Gamma_t$-module $E_t$ is the cotangent space to $T$ at $t$.
Since $\mc A_t$ is a local ring we have $\hat{\mc A_t} E_t =
\hat{\mc A_t} I_t$ by Nakayama's Lemma.
Any finite dimensional $\Gamma_t$-module is projective so
\[
U \otimes \wig^n E_t \otimes \hat{\mc A_t} =
\mr{Ind}_{\Gamma_t}^{\Gamma_t \ltimes \hat{\mc A_t}}
\big( U \otimes \wig^n E_t \big)
\]
is a projective $\Gamma_t \ltimes \hat{\mc A_t}$-module, for all
$n \in \mh N$. With these modules we construct a resolution
of $U_t$. Define $\Gamma_t \ltimes \hat{\mc A_t}$-module maps
\[
\begin{array}{lll}
\delta_n : U \otimes \wig^n E_t \otimes \hat{\mc A_t} &
\to & U \otimes \wig^{n-1} E_t \otimes \hat{\mc A_t} \\
\delta_n \, (u \otimes e_1 \wedge \cdots \wedge e_n \otimes f) & = &
\sum\limits_{i=1}^n (-1)^{i-1} u \otimes e_1 \wedge \cdots \wedge
e_{i-1} \wedge e_{i+1} \wedge \cdots \wedge e_{j_n} \otimes e_i f \\
\delta_0 : U \otimes \hat{\mc A_t} & \to & U_t \\
\delta_0 \, (u \otimes f) & = & f(t) u
\end{array}
\]
This makes
\begin{equation}\label{eq:3.9}
\big( U \otimes \wig^* E_t \otimes \hat{\mc A_t}
,\delta_* \big)
\end{equation}
into an augmented differential complex. Notice that in Mod$
\big( \hat{\mc A_t} \big)$ this just the Koszul resolution of
\[
U_t \otimes \hat{\mc A_t} \big/ I_t \hat{\mc A_t} = U_t
\]
So \eqref{eq:3.9}
is the required projective resolution of $U_t$ and
\begin{align*}
EP_{\Gamma \ltimes X} (\mr{Ind}_t & U ,\mr{Ind}_t V) \;=\;
\sum_{n=0}^\infty (-1)^n
\dim \mr{Ext}^n_{\Gamma_t \ltimes \hat{\mc A_t}} (U_t ,V_t) \\
&=\; \sum_{n=0}^r (-1)^n \dim H^n \Big( \mr{Hom}_{\Gamma_t
\ltimes \hat{\mc A_t}} \big( U \otimes \wig^* E_t
\otimes \hat{\mc A_t} , V_t \big), \mr{Hom}(\delta_* ,
\mr{id}_{V_t}) \Big) \\
&=\; \sum_{n=0}^r (-1)^n \dim \mr{Hom}_{\Gamma_t \ltimes
\hat{\mc A_t}} \big( U \otimes \wig^n E_t \otimes
\hat{\mc A_t} , V_t \big) \\
&=\; \sum_{n=0}^r (-1)^n \dim \mr{Hom}_{\Gamma_t}
\big( U \otimes \wig^n E_t , V \big) \quad
= \quad e_{\Gamma_t} (U,V)
\end{align*}
This proves c). According to Theorem \ref{thm:3.8}
$e_{\Gamma_t}$ induces an inner product on $Ell (\Gamma_t)$
and by definition $\mr{Ind}_t (Ell (\Gamma_t )) \subset Ell
(\Gamma \ltimes X)$ is precisely of the span of
the $\Gamma \ltimes X$-modules with central character
$\Gamma t$. Two $\Gamma \ltimes X$-representations
with different $Z(\Gamma \ltimes \mc A )$-characters are
orthogonal for $EP_{\Gamma \ltimes X}$, so b) and d) follow.

Now let us count the elliptic conjugacy classes in
$\Gamma \ltimes X$. Two sets\\ $(x + (\mr{id}_E - \rho
(\gamma )) X ) \gamma$ and $(y + (\mr{id}_E - \rho (\gamma )) X
) \gamma$ are conjugate if and only if there is a
$w \in Z_{\Gamma}(\gamma)$ such that $\rho (w) x - y \in
(\mr{id}_E - \rho (\gamma ) X)$. As $\Gamma$-sets we have\\
$T^\gamma = \mr{Hom} \big( X / (\mr{id}_E - \rho (\gamma )) X ,
\mh C^\times \big)$. Therefore
\begin{align*}
\# \big( (\Gamma \ltimes X)_{ell} / \sim_{\Gamma \ltimes X} \big)
&= \sum_{\gamma \in \Gamma_{ell} / \sim_\Gamma} \#
\big( (X / (1-\gamma) X ) / Z_\Gamma (\gamma) \big) \\
&= \sum_{\gamma \in \Gamma_{ell} / \sim_\Gamma} \# \big( T^\gamma /
Z_\Gamma (\gamma) \big) \\
&= \# \big( \{ (\gamma,t) : \gamma \in \Gamma_{ell} ,
t \in T^\gamma \} / Z_\Gamma (\gamma) \big) \\
&= \# \big( \{ (\gamma,t) : t \in T , \gamma \in
\Gamma_{t,ell} \} / Z_\Gamma (\gamma) \big) \\
&= \sum_{t \in T / \Gamma} \# \big( \Gamma_{t,ell} /
\sim_{\Gamma_t} \big) \\
&= \sum_{t \in T / \Gamma} \dim Ell (\Gamma_t) \\
&= \dim Ell (\Gamma \ltimes X)
\end{align*}
where we let $\Gamma$ act on $\Gamma_{ell} \times T$ by
$w \cdot (\gamma,t) = (w \gamma w^{-1} , w t). \qquad \Box$
\\[3mm]

From the above proof we see that part c) of Theorem \ref{thm:3.9}
remains valid in the following more general settings:
\begin{itemize}
\item $T$ is a nonsingular complex affine variety, $\mc A =
\mc O (T)$ and $\Gamma$ acts on $T$ by algebraic automorphisms
\item $T$ is a smooth manifold, $\mc A = C^\infty (T)$ and
$\Gamma$ acts on $T$ by diffeomorphisms.
\end{itemize}

\section{The elliptic measure}

It is shown in \cite[Theorem III.4.21]{ScSt} and
\cite[Theorem 0.20]{Bez} that the
Euler-Poincar\'e characteristic for semisimple $p$-adic groups
agrees with the elliptic integral introduced in \cite[p. 5]{Kaz}.

For the group $\Gamma \ltimes X$ this relation can be made
even more explicit. We endow it with the $\sigma$-algebra $\mc L$
generated by the sets
\begin{equation}
L_w := \{ x w (-x) : x \in X \} \qquad
w \in \Gamma \ltimes X
\end{equation}
Let $\chi_V$ denote the character of a representation $V$.

\begin{thm}\label{thm:3.10}
\begin{description}
\item[a)] There exists a unique conjugation-invariant "elliptic"
measure $\mu_{ell}$ on $(\Gamma \ltimes X ,\mc L )$ such that
\[
EP_{\Gamma \ltimes X} (U,V) = \int_{\Gamma \ltimes X}
\overline{\chi_U} \chi_V d \mu_{ell}
\qquad \forall \, U , V \in \mr{Mod}_{fin} (\Gamma \ltimes X)
\]
\item[b)] The support of $\mu_{ell}$ is the set of elliptic
elements
\item[c)] Let $e \in E$ be an isolated fixpoint of an elliptic
element $c \in \Gamma \ltimes X$ and let $C \subset \Gamma
\ltimes X$ be the conjugacy class of $c$. Then
\end{description}
\[
\begin{array}{lll}
\mu_{ell} (L_c ) & = & |\Gamma |^{-1} \\
\mu_{ell}(C) & = & {\ds \frac{\# \{ w \in C : \rho (w) e = e \} }{
\# \{ w \in \Gamma \ltimes X : \rho (w) e = e \} }} \\
\mu_{ell}(\Gamma \ltimes X ) & = & \sum\limits_{n=0}^\infty
(-1 )^n \dim \big( \wig^n E \big)^\Gamma
\end{array}
\]
\end{thm}
\emph{Proof.}
Suppose we have a trace function $f \in G_{\mh C}(\Gamma \ltimes X)$
such that $f(w) = 0 \; \forall w \in (\Gamma \ltimes X)_{ell}$.
Write $f = \sum_{t \in T / \Gamma} \mr{Ind}_t f_t$. This is a finite
sum because $G_{\mh C}(\Gamma \ltimes X)$ is built from finite
dimensional representations. If $\gamma \in \Gamma_{t,ell}$ then we
have $f (x \gamma ) = 0 \; \forall x \in X$, so
$[\Gamma : \Gamma_t ] f_t (\gamma ) = \mr{Ind}_t (f_t ) (\gamma ) = 0$.

Hence by Theorem \ref{thm:3.8}.b) $[f_t ] = 0 \in Ell (\Gamma_t )$.
By Theorem \ref{thm:3.9}.d) $[f] = 0 \in Ell (\Gamma \ltimes X)$.
Now parts and a) and
b) follow automatically, since there are only finitely many elliptic
conjugacy classes in $\Gamma \ltimes X$ and every conjugacy class
contains only finitely many $L_w$'s.

To find the explicit form of $\mu_{ell}$ we consider a possibly
different measure $\mu$ on $\Gamma \ltimes X$ defined by
$\mu (L_c ) := |\Gamma |^{-1}$,
for any elliptic element $c \in \Gamma \ltimes X$. We will show that
$\mu$ satisfies the properties attributed to $\mu_{ell}$. It will
follow from the just proven uniqueness that $\mu = \mu_{ell}$.

Let $U$ and $V$ be irreducible $\Gamma \ltimes X$-representations
with central characters $\Gamma t$ and $\Gamma t'$ respectively.
By \eqref{eq:3.2} there are characters $\chi$ of $\Gamma_t$ and
$\chi_{t'}$ of $\Gamma_{t'}$ such that $\chi_U = \mr{Ind}_t \chi$ and
$\chi_V = \mr{Ind}_t \chi'$. Extend $\chi$ and $\chi'$ to functions
on $\Gamma$ by making them zero on $\Gamma \setminus \Gamma_t$ and
on $\Gamma \setminus \Gamma_{t'}$ respectively. For $\gamma \in
\Gamma_{ell}$ we have
\[
\chi_U (x \gamma ) = \sum_{h \in \Gamma / \Gamma_t}
t (\rho (h)^{-1} x) \chi (h^{-1} \gamma h )
\]
This can only be nonzero if $\chi (h^{-1} \gamma h) \neq 0$, which
forces $h^{-1} \gamma h$ to be an elliptic element of $\Gamma_t$.
Therefore
\[
\int_{\Gamma \ltimes X} \overline{\chi_U} \chi_V d \mu = 0
\]
if either $Ell (\Gamma_t ) = 0$ or $Ell (\Gamma_{t'}) = 0$, which is
in agreement with Theorem \ref{thm:3.8}.b).

Hence we assume that $\Gamma_t$ and $\Gamma_{t'}$ do contain elliptic
elements. This forces all elements of $\Gamma \{t , t'\}$ to have
finite order in the group $T$. Now
\[
X' := \bigcap_{t'' \in \Gamma \{t , t'\}} \ker t'' \cap
\bigcap_{\gamma \in \Gamma_{ell}} (\mr{id}_E - \rho (\gamma ) )X
\]
is a lattice of finite index in $X$ and the map
\[
X / X' \to \mh C : x \to t( \rho (h)^{-1} x) \chi (h^{-1} \gamma h)
\]
is well defined for all $h,\gamma \in \Gamma$. For a fixed
$\gamma \in \Gamma_{ell}$ we have
\begin{equation}\label{eq:3.13}
\begin{aligned}
{[} (\mr{id}_E - & \rho (\gamma ) )X : X' {]}
\!\!\!\! \sum_{x \in X / (\mr{id}_E - \rho (\gamma ) )X} \!\!\!\!
\overline{\chi_U (x \gamma )} \chi_V (x \gamma ) = \sum_{x \in X / X'}
\overline{\chi_U (x \gamma )} \chi_V (x \gamma ) \\
& = \sum_{h \in \Gamma / \Gamma_t} \sum_{g \in \Gamma / \Gamma_{t'}}
\sum_{x \in X / X'} \overline{t (\rho (h)^{-1} x) \chi (h^{-1}
\gamma h )} \: t' (\rho (g)^{-1} x) \chi' (g^{-1} \gamma g )
\end{aligned}
\end{equation}
By the orthogonality relations for characters of the group
$X / X'$ the only nonzero contributions to this sum come from pairs
$(g,h)$ for which $h (t) = g (t')$. In particular
\[
\int_{\Gamma \ltimes X} \overline{\chi_U} \chi_V d \mu = 0
\]
if $\Gamma t \neq \Gamma t'$.
This leaves the case $t = t'$. From \eqref{eq:3.13} we see that
\begin{align*}
\sum_{x \in X / (\mr{id}_E - \rho (\gamma ) )X}
\!\!\!\!\!\!\!\!\!\!\!\!\! \overline{\chi_U (x \gamma )}
\chi_V (x \gamma ) \, & = \sum_{h,g \in \Gamma / \Gamma_t}
\sum_{x \in X / X'} \!\!\! \frac{\overline{t( \rho (h)^{-1} x) \chi (h^{-1}
\gamma h )} \: t( \rho (g)^{-1} x) \chi' (g^{-1} \gamma g )}{
[ (\mr{id}_E - \rho (\gamma ) )X : X' ]} \\
& = \sum_{h \in \Gamma / \Gamma_t} \sum_{x \in X / X'} \frac{\overline{t(
\rho (h)^{-1} x) \chi (h^{-1} \gamma h )} \: t( \rho (h)^{-1} x) \chi'
(h^{-1} \gamma h )}{[ (\mr{id}_E - \rho (\gamma ) )X : X' ]} \\
& = [X : (\mr{id}_E - \rho (\gamma ) )X] \sum_{h \in \Gamma / \Gamma_t}
\overline{\chi (h^{-1} \gamma h )} \chi' (h^{-1} \gamma h )
\end{align*}
Now we can compute
\begin{align*}
\int_{\Gamma \ltimes X} \overline{\chi_U} \chi_V d \mu \, & =
\sum_{\gamma \in \Gamma_{ell}} \sum_{x \in X / (\mr{id}_E - \rho (\gamma ) )X}
\frac{\overline{\chi_U (x \gamma )} \chi_V (x \gamma )}{|\Gamma |} \\
& = \sum_{\gamma \in \Gamma_{ell}} \frac{[X :
(\mr{id}_E - \rho (\gamma ) )X]}{|\Gamma |} \sum_{h \in \Gamma / \Gamma_t}
\overline{\chi (h^{-1} \gamma h )} \chi' (h^{-1} \gamma h ) \\
& = \sum_{\gamma \in \Gamma_{t,ell}} \frac{\det \, (\mr{id}_E - \rho (\gamma ))
}{|\Gamma |} [\Gamma : \Gamma_t ] \overline{\chi (\gamma )} \chi' (\gamma ) \\
& = \sum_{\gamma \in \Gamma_{t,ell}} \frac{\det \, (\mr{id}_E - \rho (\gamma ))
}{|\Gamma_t |} \overline{\chi (\gamma )} \chi' (\gamma ) \quad
= \quad e_{\Gamma_t} (\chi ,\chi') \\
\end{align*}
Thus indeed $\mu = \mu_{ell}$.

Let $e,c$ and $C$ be as above. To determine $\mu_{ell}(C)$ we must count
the number $n_C$ of sets $L_w$ that are contained in $C$. Consider the map
\begin{align*}
& \psi_e : C \to E / X \\
& \psi_e (w c w^{-1}) = \rho(w) e + X
\end{align*}
It is easily seen that $\psi_e$ is well-defined and that
\[
\psi_e (x w c w^{-1} (-x)) = \psi_e (w c w^{-1})
\qquad \forall x \in X ,\, w \in \Gamma \ltimes X
\]
The image of $\psi_e$ is $\rho (\Gamma \ltimes X) e / X$ and
\[
\psi_e^{-1}(\rho (w) e + X ) = \{ x w v c v^{-1} w^{-1} (-x) :
x \in X , v \in \Gamma \ltimes X , \rho (v) e = e \}
\]
The number of $L_w$'s contained in $\psi_e^{-1} (\rho (w) e + X)$ is
\[
\# \{ v c v^{-1} : v \in \Gamma \ltimes X , \rho (v) e = e \} =
\# \{ v \in C : \rho (v) e = e \}
\]
Consequently
\begin{align*}
& n_C = |\rho (\Gamma \ltimes X) e / X| \: \# \{ v \in C : \rho (v) e = e \}
= \frac{|\Gamma | \: \# \{ v \in C : \rho (v) e = e \}}{\#
\{ w \in \Gamma \ltimes X : \rho (w) e = e \}} \\
& \mu_{ell}(C) = \frac{n_C}{|\Gamma |} = \frac{\# \{ v \in C :
\rho (v) e = e \}}{\# \{ w \in \Gamma \ltimes X : \rho (w) e = e \}}
\end{align*}
Finally, using Theorem \ref{thm:3.9}.c) we compute
\begin{align*}
\mu_{ell}(\Gamma \ltimes X) & = EP_{\Gamma \ltimes X}
(\mr{triv}_{\Gamma \ltimes X} ,\mr{triv}_{\Gamma \ltimes X} ) \\
& = EP_{\Gamma \ltimes X} \big( \mr{Ind}_1 (\mr{triv}_\Gamma) ,
\mr{Ind}_1 (\mr{triv}_\Gamma) \big) \\
& = e_\Gamma \big( \mr{triv}_\Gamma, \mr{triv}_\Gamma \big) \\
& = \sum\limits_{n=0}^\infty (-1)^n \dim \mr{Hom}_\Gamma \big(
\wig^n E , \mr{triv}_\Gamma \big) \\
& = \sum\limits_{n=0}^\infty (-1)^n
\dim \big( \wig^n E \big)^\Gamma \qquad \Box
\end{align*}
\\[1mm]

\section{Example: the Weyl group of type $B_2$}

Let $R_0$ be the root system $B_2$ in $E = \mh R^2$,
with positive roots
\[
\alpha_1 = (1,-1) ,\, \alpha_2 = (0,1) ,\, \alpha_3 (1,0)
,\, \alpha_4 = (1,1)
\]
Denote the rotation of $E$ over an angle $\theta$ by $\rho_\theta$
and the reflection corresponding to $\alpha_i$ by $s_i$. Then
\[
W_0 = \{ e, s_1 ,s_2 ,s_3 ,s_4 , \rho_{\pi /2}, \rho_\pi,
\rho_{-\pi /2} \}
\]
is isomorphic to the dihedral group $D_4$. This group has four
irreducible representations of dimension one, defined by
\begin{equation}
\begin{array}{c|rr}
\pi & \pi (s_1 ) & \pi (s_2) \\
\hline
\ep_0 & 1 & 1 \\
\ep_1 & -1 & 1 \\
\ep_2 & 1 & -1 \\
\ep_3 & -1 & -1
\end{array}
\end{equation}
The one remaining irreducible representation is just $E$.

The elliptic conjugacy classes in $W_0$ are $\{ \rho_\pi \}$ and
$\{ \rho_{\pi / 2} , \rho_{-\pi / 2} \}$.
\[
\begin{array}{lll}
\mr{Ind}_{W_\es}^{W_0} (G_{\mh C} \{e\}) & = & \mh C \{
\ep_0 \oplus \ep_1 \oplus \ep_2 \oplus \ep_3 \oplus E \oplus E \} \\
\mr{Ind}_{W_{\{1\}}}^{W_0} (G_{\mh C} \{e,s_1 \}) & = & \mh C \{
\ep_0 \oplus \ep_2 \oplus E \,, \ep_1 \oplus \ep_3 \oplus E \} \\
\mr{Ind}_{W_{\{2\}}}^{W_0} (G_{\mh C} \{e,s_2 \}) & = & \mh C \{
\ep_0 \oplus \ep_1 \oplus E \,, \ep_2 \oplus \ep_3 \oplus E \}
\end{array}
\]
We see that $Ell (W_0 )$ has dimension two and is spanned for
example by $[\ep_0 ]$ and $[\ep_1 ]$. With Theorem \ref{thm:3.8}.c)
we can easily write down a complete table for $e_{W_0}$ :
\begin{equation}
\begin{array}{c|rrrrr}
e_{W_0} & \ep_0 & \ep_1 & \ep_2 & \ep_3 & E \\
\hline
\ep_0 & 1 & 0 & 0 & 1 & -1 \\
\ep_1 & 0 & 1 & 1 & 0 & -1 \\
\ep_2 & 0 & 1 & 1 & 0 & -1 \\
\ep_3 & 1 & 0 & 0 & 1 & -1 \\
E & -1 & -1 & -1 & -1 & 2
\end{array}
\end{equation}
Since $\overline{A_\es}$ is a fundamental domain for the action
of $W$ on $E$, every point of $E$ that is fixed by an elliptic
element of $W$ must be in the $W$-orbit of some vertex of the
fundamental alcove $A_\es$. This leads to the following list
of elliptic conjugacy classes:
\begin{equation}
\begin{array}{ccc}
\text{vertex} & \text{conjugacy class} & \text{elliptic measure} \\
e = c (e) & [c] & \mu_{ell}([c]) \\
\hline
(0,0) & [\rho_\pi ] & 1/8 \\
(0,0) & [\rho_{\pi /2}] & 1/4 \\
(1/2,1/2) & [t_{(1,1)} \rho_\pi ] & 1/8 \\
(1/2,1/2) & [t_{(1,0)} \rho_{\pi /2} ] & 1/4 \\
(1/2,0) & [t_{(1,0)} \rho_\pi ] & 1/4
\end{array}
\end{equation}
In particular dim $Ell (W) = 5$.

For $t \in T$ we write $t = (t (1,0) , t (0,1) )$. The following
points of $T$ are fixed by an elliptic element of $W_0$.
\begin{itemize}
\item $(1,1)$ is fixed by all $w \in W_0$. Thus we get a
two dimensional subspace $\mr{Ind}_{(1,1)} \big( Ell (W_0 ) \big)$
of $Ell (W)$.
\item $(-1,-1)$ is also fixed by the whole group $W_0$. This gives
another two dimensional subspace $\mr{Ind}_{(-1,-1)} \big(
Ell (W_0 ) \big) \subset Ell (W)$.
\item $(-1,1)$ has isotropy group $V_4 = \{ e,s_2 ,s_3 ,\rho_\pi \}
\subset W_0$. The only elliptic element is $\rho_\pi$ so dim
$Ell (V_4 ) = 1$.
\item $(1,-1)$ also has isotropy group $V_4$. But $(-1,1)$ and
$(1,-1)$ are in the same $W_0$-orbit so
$\mr{Ind}_{(1,-1)} \big( Ell (V_4 ) \big) =
\mr{Ind}_{(-1,1)} \big( Ell (V_4 ) \big)$.
This one dimensional subspace of $Ell (W)$ is spanned for example
by the two dimensional representation
$\mr{Ind}_{(1,-1)}\big( \mr{triv}_{V_4 } \big)$.
\end{itemize}
Now we have three subspaces of $Ell (W)$ that are mutually
orthogonal for $EP_W$ and whose dimensions add up to 5. Since
this is exactly the number of elliptic conjugacy classes in $W$,
we found all of $Ell (W)$.
\\[2mm]

\section{The Euler-Poincar\'e characteristic}

Following Schneider and Stuhler \cite[\S III.4]{ScSt} we introduce
an Euler-Poincar\'e characteristic for affine Hecke algebras.
For finite dimensional $\mc H$-modules $U$ and $V$ we define
\begin{equation}
EP_{\mc H}(U,V) = \sum_{n=0}^\infty (-1)^n \dim
\mr{Ext}^n_{\mc H}(U,V)
\end{equation}
By Proposition \ref{prop:2.dimh} the sum is actually finite,
so this is well-defined.
With standard homological algebra (see for instance \cite{CaEi} ) one
can show that this extends to a bilinear pairing on $G (\mc H)$.
Reeder \cite{Ree} studied this pairing for affine Hecke algebras
with equal parameters, via $p$-adic groups.

\begin{prop}\label{prop:3.6}
\begin{description}
\item[a)] Let $I \subset F_0$ be a proper subset of simple roots and let\\
$V \in \mr{Mod}_{fin} \big( \mc H (\mc R^I ,q^I) \big)$. Then
\[
EP_{\mc H} \big( U, \mr{Ind}_{\mc H (\mc R^I ,q^I )}^{\mc H} V \big) = 0
\qquad \forall \, U \in \mr{Mod}_{fin}(\mc H )
\]
\item[b)] If the root datum $\mc R$ is not semisimple then
$EP_{\mc H} \equiv 0$.
\end{description}
\end{prop}
\emph{Proof.}
This result is the translation of \cite[Lemma III.4.18.ii]{ScSt} to
affine Hecke algebras. The proof is similar and based on an argument
due to Kazhdan.

We may assume that $(\pi ,V)$ is irreducible with $Z(\mc H (\mc R^I ,q^I)
)$-character $W_I t \in T / W_I$. If $W_0 t$ is not an $Z (\mc H )$-weight
of $U$ then $\mr{Ext}_n^{\mc H} \big( U, \mr{Ind}_{\mc H (\mc R^I ,q^I
)}^{\mc H} V \big) = 0$, so certainly $EP_{\mc H}\big( U,\mr{Ind}_{\mc H
(\mc R^I ,q^I )}^{\mc H} V \big) = 0$. Therefore we may also assume that
$U$ is irreducible with $Z (\mc H )$-character $W_0 t \in T /W_0$.
Recall the group $\tilde G$ from \eqref{eq:2.29} and abbreviate
\[
\begin{array}{lll}
m_t & = & \# \{ g \in \widetilde G : g W_0 t = W_0 t \} \\
\widetilde U & = & \mr{Ind}_{\mc H (\mc R ,q)}^{
\mc H (\widetilde{\mc R},\tilde q)} (U)
\end{array}
\]
Let $\mc F_n$ be a set of representatives for the action of
$W_0 (\widetilde{\mc R})$ on the facets of the fundamental alcove for
$\widetilde{\mc R}$.
From \eqref{eq:2.28} and Corollary \ref{cor:2.11}.a) we deduce that
\begin{align}
\nonumber m_t \: & EP_{\mc H} \big( U ,
\mr{Ind}_{\mc H (\mc R^I ,q^I )}^{\mc H} V \big) \; = \;
EP_{\mc H} \big( U , \mr{Ind}_{\mc H (\mc R^I ,q^I )}^{
\mc H (\widetilde{\mc R},\tilde q)} V \big) \\
\nonumber & = \; EP_{\mc H (\widetilde{\mc R},\tilde q)} \big( \widetilde U ,
\mr{Ind}_{\mc H (\mc R^I ,q^I )}^{\mc H (\widetilde{\mc R},\tilde q)} V
\big) \\
\label{eq:3.1} & = \; \sum_{n=0}^\infty (-1)^n \dim \mr{Ext}^n_{\mc H
(\widetilde{\mc R} ,\tilde q)} \big( \widetilde U , \mr{Ind}_{\mc H
(\mc R^I ,q^I )}^{\mc H (\widetilde{\mc R} ,\tilde q)} V \big) \\
\nonumber & = \; \sum_{n=0}^{\mr{rk}(X)} (-1)^n \dim H^n \Big(
\bigoplus_{f \in \mc F_n} \mr{Hom}_{\mc H (\widetilde{\mc R},f
,\tilde q)} \big( \widetilde U \otimes \ep_f , \mr{Ind}_{\mc H (\mc R^I
,q^I )}^{\mc H (\widetilde{\mc R} ,\tilde q)} V \big)
,\mr{Hom}(d_* ,\mr{id} ) \Big) \\
\nonumber & \; = \sum_{n=0}^{\mr{rk}(X)} (-1)^n \sum_{f \in \mc F_n}
\dim \mr{Hom}_{\mc H (\widetilde{\mc R},f ,\tilde q)}
\big( \widetilde U \otimes \ep_f , \mr{Ind}_{\mc H (\mc R^I ,q^I )
}^{\mc H (\widetilde{\mc R} ,\tilde q)} V \big)
\end{align}
Because $V$ is irreducible there exist a $\mc H (\mc R_I ,q_I
)$-representation $(\pi_1 ,V)$ and a\\
$Z \big( W(\mc R^I ) \big)$-character $t_1$ such that
\[
(\pi ,V) = (\pi_1 \circ \phi_{t_1} ,V)
\]
with $\phi_{t_1}$ as in \eqref{eq:2.19}. Note that
$Z \big( W(\mc R^I ) \big) = (I^\vee )^\perp \cap X \neq 0$ because
$I \neq F_0$. Let $t_2$ be an arbitrary $Z \big( W(\mc R^I )
\big)$-character and consider the integer
\[
\dim \mr{Hom}_{\mc H (\widetilde{\mc R},f ,\tilde q)} \big( \widetilde U
\otimes \ep_f , \mr{Ind}_{\mc H (\mc R^I ,q^I )}^{\mc H (\widetilde{\mc R}
,\tilde q)} (\pi_1 \circ \phi_{t_2} , V ) \big)
\]
According to Lemma \ref{lem:2.5} $\mc H (\widetilde{\mc R},f,\tilde q)$
is a finite dimensional semisimple algebra. Therefore the above integer is
invariant under continuous deformations of $t_2$, and hence independent
of $t_2$. Pick $t_2$ such that the central character of \\
$\mr{Ind}_{\mc H (\mc R^I ,q^I )}^{\mc H (\widetilde{\mc R} ,\tilde q)}
(\pi_1 \circ \phi_{t_2} , V ) \text{  is not  }
W_0 (\widetilde{\mc R}) t \in T / W_0 (\widetilde{\mc R})$. Then
\begin{equation}
\begin{split}
0 & = m_t EP_{\mc H} \big( U , \mr{Ind}_{\mc H (\mc R^I ,q^I )}^{\mc H
(\widetilde{\mc R} ,\tilde q)} (\pi_1 \circ \phi_{t_2} , V ) \big) \\
& = \sum_{n=0}^{\mr{rk}(X)} (-1)^n \sum_{f \in \mc F_n} \dim
\mr{Hom}_{\mc H (\widetilde{\mc R},f ,\tilde q)} \big( \widetilde U
\otimes \ep_f , \mr{Ind}_{\mc H (\mc R^I ,q^I )}^{\mc H
(\widetilde{\mc R} ,\tilde q)} (\pi_1 \circ \phi_{t_2} , V ) \big) \\
& = \sum_{n=0}^{\mr{rk}(X)} (-1)^n \sum_{f \in \mc F_n} \dim
\mr{Hom}_{\mc H (\widetilde{\mc R},f ,\tilde q)} \big( \widetilde U
\otimes \ep_f , \mr{Ind}_{\mc H (\mc R^I ,q^I )}^{\mc H
(\widetilde{\mc R} ,\tilde q)} (\pi_1 \circ \phi_{t_1} , V ) \big) \\
& = m_t EP_{\mc H} \big( U , \mr{Ind}_{\mc H (\mc R^I ,q^I )}^{\mc H}
(\pi , V) \big)
\end{split}
\end{equation}
To prove b) we suppose that $\mc R$ is not semisimple and that
$U' , V' \in \mr{Mod}_{fin}(\mc H )$. We have to show that
\[
EP_{\mc H} (U' ,V' ) = 0
\]
We may assume that $U'$ and $V'$ admit the same central character $W_0 t$.
From the proof of part a) we see that
\[
m_t \, EP_{\mc H} (U' ,V' ) = EP_{\mc H (\widetilde{\mc R},\tilde q)}
\big( \mr{Ind}_{\mc H}^{\mc H (\widetilde{\mc R} ,\tilde q)} U' ,
\mr{Ind}_{\mc H}^{\mc H (\widetilde{\mc R},\tilde q)} V' \big) = 0
\qquad \Box
\]
\vspace{2mm}

We can use the scaling maps
\[
\tilde \sigma_\ep : \mr{Mod}_{fin}(\mc H (\mc R ,q)) \to
\mr{Mod}_{fin}(\mc H (\mc R ,q^\ep ))
\]
from Theorem \ref{thm:2.4} to relate $EP_{\mc H}$ to $EP_W$.

\begin{thm}\label{thm:3.4}
\begin{description}
\item[a)] The pairing $EP_{\mc H}$ is symmetric and positive semidefinite.
\item[b)] If $U , V \in \mr{Mod}_{fin}(\mc H )$ then
\[
EP_{\mc H}(U,V) = EP_{\mc H (\mc R ,q^\ep )} \big( \tilde \sigma_\ep (U),
\tilde \sigma_\ep (V) \big) \qquad \forall \ep \in [-1,1]
\]
\end{description}
\end{thm}
\emph{Proof.}
In view of Proposition \ref{prop:3.6}.b) we may assume that $\mc R$ is
semisimple. For every $\ep \in [-1,1]$ Theorem \ref{thm:2.4} gives us the
$\mc H (\mc R,q^\ep )$-representations
\[
\tilde \sigma_\ep (\rho ,U) = (\rho_\ep ,U) \qquad \mr{and} \qquad
\tilde \sigma_\ep (\pi ,V) = (\pi_\ep ,V) .
\]
As a vector space $\mc H (\mc R ,f,q^\ep )$ is just $\mh C [W_f \rtimes
\Omega_f ]$. As an algebra it is semisimple and the multiplication
varies continuously with $\ep$, so by Tits' deformation lemma it is
independent of $\ep$. Furthermore for any $w \in W_f \rtimes \Omega_f$
the maps
\[
\ep \to \rho_\ep (N_w ) \qquad \mr{and} \qquad \ep \to \pi_\ep (N_w )
\]
are continuous. In view of \eqref{eq:3.1} this implies that
\[
EP_{\mc H (\mc R ,q^\ep )} \big( \tilde \sigma_\ep (U),\tilde \sigma_\ep
(V) \big)
\]
depends continuously on $\ep$. But this expression is integer valued,
so it is actually independent of $\ep$. In particular
\begin{equation}\label{eq:3.EP}
EP_{\mc H} (U,V) = EP_W \big( \tilde \sigma_0 (U),\tilde \sigma_0 (V) \big)
\end{equation}
Now Theorem \ref{thm:3.9}.b) assures that $EP_{\mc H}$ is
symmetric and positive semidefinite. $\qquad \Box$
\\[3mm]

For semisimple root data we can also compute the Euler-Poincar\'e
characteristic in another way, as the character value of
a certain index function.

According to Lemma \ref{lem:2.5} for all facets $f$ the algebra
$\mc H (\mc R ,f,q)$ is finite dimensional and semisimple, so in
particular the collection Irr$ (\mc H (\mc R ,f,q))$ of irreducible
representations is finite. Let $e_\sigma \in \mc H (\mc R ,f,q)$
denote the primitive central idempotent corresponding to an
irreducible $\mc H (\mc R ,f,q)$-module $\sigma$.
For $U \in \mr{Mod}(\mc H (\mc R ,f,q))$ let $[U : \sigma ]$ be the
multiplicity of $\sigma$ in $U$.

In the spirit of Kottwitz \cite[\S 2]{Kot}, Schneider and Stuhler
\cite[III.4]{ScSt} we define an Euler-Poincar\'e function
\begin{equation}
f_{EP}^U := \sum_{f \subset \overline{A_\es}} \frac{(-1 )^{\dim f}
}{[\Omega : \Omega_f ]} \sum_{\sigma \in \mr{Irr} (\mc H (\mc R ,f,q) )}
\frac{[U \otimes \ep_f : \sigma]}{\dim \sigma} e_\sigma
\end{equation}

\begin{prop}\label{prop:3.3}
Let $\mc R$ be a semisimple root datum and
$U,V \in \mr{Mod}_{fin}(\mc H )$. Then
\[
EP_{\mc H} (U,V) = \chi_V \big( f_{EP}^U \big)
\]
\end{prop}
\emph{Proof.}
Exactly like in \eqref{eq:3.1} we can calculate that
\begin{align*}
EP_{\mc H} (U,V) & = \sum_{n = 0}^r (-1 )^n \dim \mr{Hom}_{\mc H}
(P_n (U )^\Omega ,V) \\
& = \sum_{n = 0}^r \sum_{f : \dim f = n} \frac{(-1 )^n
}{[\Omega : \Omega_f ]} \dim \mr{Hom}_{\mc H (\mc R ,f,q)}
(U \otimes \ep_f , V ) \\
& = \sum_{f \subset \overline{A_\es}}  \frac{(-1 )^{\dim f}}{[
\Omega : \Omega_f ]} \sum_{\sigma \in \mr{Irr} (\mc H (\mc R ,f,q) )}
[U \otimes \ep_f : \sigma ] \, [V : \sigma ] \\
& = \sum_{f \subset \overline{A_\es}}  \frac{(-1 )^{\dim f}}{[
\Omega : \Omega_f ]} \sum_{\sigma \in \mr{Irr} (\mc H (\mc R ,f,q) )}
\frac{[U \otimes \ep_f : \sigma]}{\dim \sigma} \chi_V (e_\sigma ) \\
& = \chi_V \big( f_{EP}^U \big) \hspace{5cm} \Box
\end{align*}
\\[2mm]
We will use this result in \cite{OpSo} to show that the Plancherel
measure of a discrete series representation is a rational function
in $q$, with rational coefficients.
\vspace{4mm}

\section{Extensions of tempered modules}

We apply the results of Chapter 2 to relate the
bornological Tor and Ext functors over $\mc H$ with
those over $\mc S$. That is more interesting than it looks at
first sight, because $\mc S$ is not flat over $\mc H$ (unless
$q \equiv 1$).

\begin{cor}\label{cor:3.2}
Take $n \in \mh N$.
\begin{description}
\item[a)] For all $U_b ,V_b \in \mr{Mod}_{bor} (\mc S )$
the inclusion $\mc H \to \mc S$ induces isomorphisms
\[
\begin{array}{lllll}
\mr{Tor}_n^{\mc H} (\mc S ,V_b) & \cong &
\mr{Tor}_n^{\mc S} (\mc S ,V_b) & \cong &
\left\{
\begin{array}{lll}
V_b & \mr{if} & n = 0 \\
0 & \mr{if} & n > 0
\end{array} \right. \\
\mr{Ext}^n_{\mc H}(U_b ,V_b ) & \cong & \mr{Ext}^n_{\mc S}(U_b ,V_b )
\end{array}
\]
\item[b)] For all finite dimensional tempered $\mc H$-modules
$U$ and $V$ there is a natural\\ isomorphism
$\mr{Ext}^n_{\mc H}(U,V) \cong \mr{Ext}^n_{\mc S}(U,V) \,$.
\item[c)] $EP_{\mc H} (U,V) = EP_{\mc S} (U,V) \,$.
\end{description}
\end{cor}
\emph{Proof.}
a) follows directly from Theorems \ref{thm:2.16}.a) and \ref{thm:2.15}.\\
b) In this setting the bornological functor $\mr{Ext}_n^{\mc H}$
agrees with its purely algebraic counterpart, as discussed in the
appendix. The same holds for $\mr{Ext}_n^{\mc S}$, because the
resolution from Corollary \ref{cor:2.9} consists of $\mc S$-modules
that are projective in both the algebraic and the bornological sense.
Hence b) is a special case of a).

However for semisimple root data this can be proved more directly,
without the use of bornological techniques. Namely, we can simply
compare the projective resolutions from Corollaries \ref{cor:2.11}.a)
and \ref{cor:2.9}.a). If we use these to compute the Ext-groups and
we apply Frobenius reciprocity, then we see that
$\mr{Ext}^n_{\mc H}(U,V)$ and $\mr{Ext}^n_{\mc S}(U,V)$ are the
homologies of isomorphic differential complexes.
See also the remark at the end of Section \ref{sec:2.3}.\\
c) is a trivial consequence of b). $\qquad \Box$
\\[3mm]

Notice that we have to take the derived functors with respect
to bornological tensor products and bounded maps if we want to
get Corollary \ref{cor:3.2}.a) for infinite dimensional modules.
If we would work purely algebraically this would already fail
for $U = V = \mc S$.

The main use of Corollary \ref{cor:3.2}.c) is the next theorem.
Notice that the proof of the corresponding result for reductive
$p$-adic groups \cite[Theorem 41]{Mey2} is a lot more involved.

\begin{thm}\label{thm:3.5}
Suppose that $U$ and $V$ are irreducible tempered $\mc H$-modules.
If $U$ or $V$ belongs to the discrete series then
\[
\mr{Ext}_{\mc H}^n (U,V) \cong \left\{ \begin{array}{lll}
\mh C & \mr{if} \; U \cong V \; \mr{and} \; n = 0 \\
0 & \mr{otherwise}
\end{array} \right.
\]
\end{thm}
\emph{Proof.}
The assertion for $n = 0$ follows directly from Schur's lemma and the
general isomorphism $\mr{Ext}^0 \cong \mr{Hom}$.

Let $\delta$ be a discrete series representation of $\mc H$. According
to \cite[Corollary 3.13]{DeOp} $\mr{End}_{\mh C}(\delta )$ is a direct
summand of $\mc S$, as algebras. Therefore $\delta$ is both injective
and projective as a $\mc S$-module. Thus for any tempered
$\mc H$-module $V$ and any $n > 0$ we have
\begin{equation}
\mr{Ext}_{\mc H}^n (V,\delta ) = \mr{Ext}_{\mc S}^n (V,\delta ) = 0
\end{equation}
because $\delta$ is injective and
\begin{equation}
\mr{Ext}_{\mc H}^n (\delta ,V) = \mr{Ext}_{\mc S}^n (\delta ,V) = 0
\end{equation}
because $\delta$ is projective. $\qquad \Box$
\\[3mm]

Let us introduce the space of "elliptic trace functions"
\begin{equation}\label{eq:3.11}
Ell (\mc H ) := G_{\mh C}(\mc H) \Big/ \sum_{I \subset F_0 ,
I^\perp \neq 0} \mr{Ind}_{\mc H (\mc R^I ,q^I )}^{\mc H} G_{\mh C}
\big( \mc H (\mc R^I ,q^I ) \big)
\end{equation}
where $I^\perp = \{ y \in Y : \inp{\alpha}{y} = 0 \: \forall \alpha \in
I \}$. Notice that this space is zero whenever $\mc R$ is not semisimple.
From Proposition \ref{prop:3.6} and Theorem \ref{thm:3.4} we see that
the Euler-Poincar\'e characteristic induces a semidefinite Hermitian form
on $Ell (\mc H)$:
\[
EP_{\mc H} (\lambda [U] ,\mu [V]) := \bar \lambda \mu \, EP_{\mc H} (U,V)
\qquad \qquad U,V \in \mr{Mod}_{fin}(\mc H) ,\, \lambda, \mu \in \mh C
\]

\begin{prop}\label{prop:3.9}
\begin{description}
\item[a)] The scaling map $\tilde \sigma_0$ induces a linear map
$Ell (\mc H ) \to Ell (W)$ which is an isometry with respect
to the (semidefinite) Hermitian forms $EP_{\mc H}$ and $EP_W$.
\item[b)] The number of inequivalent discrete series
representations of $\mc H$ is at most the number of elliptic
conjugacy classes in $W$.
\end{description}
\end{prop}
\emph{Proof.}
a) follows directly from Theorem \ref{thm:3.4}.b)\\
b) According to Theorem \ref{thm:3.5} the inequivalent discrete series
representations form an orthonormal set in $Ell (\mc H )$. By part a)
the same holds for their images in $Ell (W)$. From
Theorem \ref{thm:3.9}.a) we know that the dimension of $Ell (W)$
is precisely the number of elliptic conjugacy classes in $W .
\qquad \Box$
\\[3mm]
\textbf{Remark.}\\
A lower bound for the number of discrete series representations can be
obtained from counting their central characters. In turns out that for the
crucial irreducible non-simply laced cases $C_n^{(1)}, F_4$ and $G_2$
this lower bound equals the above upper bound, for \emph{generic} parameters.
We will exploit this in \cite{OpSo} to give a classification of the
irreducible discrete series characters for any irreducible non-simply laced
affine Hecke algebra, with \emph{arbitrary} positive parameters.
\vspace{2mm}

\noindent\textbf{Example.}\\
Let $R_0 = A_1 = \{ 1,-1 \}$ and $X = \mh Z$. Then $W_0 = \{e,s\}$ and
$W$ is generated by $s$ and $t_1 s$. Take a label function such that
$q(s) = q (t_1 s) = q > 1$. The affine Hecke algebra $\mc H (A_1 ,q)$
has a unique discrete series representation called the Steinberg
representation. It has dimension one and is defined simply by
\[
\mr{St} (N_w) = (-1)^{\ell (w)} q(w)^{-1/2} = (-q^{-1/2})^{\ell (w)}
\]
On the other hand we have the "trivial" $\mc H$-representation
defined by
\[
\mr{triv}_{\mc H} (N_w ) = q(w)^{1/2} = q^{\ell (w) /2}
\]
It is unitary but not tempered. From Theorem \ref{thm:3.5}
we see that
\[
EP_{\mc H}(\mr{St},\mr{St}) = 1
\]
but is not immediately clear how many extensions of $\mr{St}$ by
$\mr{triv}_{\mc H}$ there are. There certainly is an extension
\begin{equation}\label{eq:3.10}
0 \leftarrow \mr{St} \leftarrow \mr{Ind}_{\mc A}^{\mc H}
(\phi_{q^{-1}}) \leftarrow \mr{triv}_{\mc H} \leftarrow 0
\end{equation}
so $\big[ \mr{Ind}_{\mc A}^{\mc H} (\phi_{q^{-1}}) \big] =
[\mr{St}] + [\mr{triv}_{\mc H}]$ in $G (\mc H )$. Therefore
\begin{align*}
EP_{\mc H} (\mr{St},\mr{triv}_{\mc H}) &= EP_{\mc H} \big( \mr{St},
[\mr{triv}_{\mc H}] - \big[ \mr{Ind}_{\mc A}^{\mc H} (\phi_{q^{-1}})
\big] \big) \\
&= EP_{\mc H} \big( \mr{St}, -[\mr{St}] \big) = -1
\end{align*}
From Corollary \ref{cor:2.11}.d) we know that the cohomological
dimension of Mod$ (\mc H )$ is 1, so in particular
\[
\mr{Ext}^n_{\mc H} (\mr{St},\mr{triv}_{\mc H}) = 0
\quad \mr{for} \; n > 1.
\]
Therefore \eqref{eq:3.10} is up to a scalar factor the only
nontrivial extension of $\mr{St}$ by $\mr{triv}_{\mc H}$.

%% file: temperedA.tex
\appendix
\chapter{Bornological algebras}

In the chapter 2 we induce several modules from $\mc H$
to $\mc S$. From an analytical point of view this operation is
trivial for finite dimensional modules, since in that case
all involved tensor products are purely algebraic.
However for infinite dimensional modules we have to take the
topology into account. For Fr\'echet $\mc S$-modules we can
use the complete projective tensor product. But for tensor
products over $\mc H$ this is problematic as there is no
canonical topology on $\mc H$.

Consider for example the trivial onedimensional root datum
$(\mh Z ,\es ,\mh Z ,\es )$. Then
\[
\mc H = \mh C [\mh Z ] \cong \mc O \big( \mh C^\times \big) \quad ,
\quad
\mc S = \mc S (\mh Z ) \cong C^\infty \big( S^1 \big) \,.
\]
For $t \in S^1$ the ideal
\[
J_t := \{ f \in C^\infty \big( S^1 \big) : f(t) = 0 \}
\subset C^\infty \big( S^1 \big)
\]
is generated by $J_t \cap \mc O \big( \mh C^\times \big)$. It
follows that for any finite dimensional $\mc S (\mh Z )$-module
$V$ we have
\[
\mc S (\mh Z ) \otimes_{\mh C [\mh Z ]} V \cong
\mc S (\mh Z ) \otimes_{\mc S (\mh Z )} V = V \,.
\]
This property does not readily generalize to infinite dimensional
modules, for example
\[
\mc S (\mh Z ) \otimes_{\mh C [\mh Z ]} \mc S (\mh Z ) \not\cong
\mc S (\mh Z ) \otimes_{\mc S (\mh Z )} \mc S (\mh Z ) =
\mc S (\mh Z ) \,.
\]
The right technique to fix this is bornology. On many vector spaces
bornological and topological analysis are equivalent, but
bornologies combine well with homological algebra in larger classes.
Bornologies are not so well-known, so we provide a brief
introduction. See also \cite{Mey0,Mey1}.

A bornology on a complex vector space is a certain collection
of subsets
that are called bounded. This collection has to satisfy some
axioms that generalize obvious properties of bounded sets in
Banach spaces. A morphism of bornological vector spaces is a
linear map that sends bounded sets to bounded sets. There is
a natural notion of completeness of bornological vector spaces,
similar to that of completeness of locally convex spaces.

On any vector space $V$ we can define a more or less trivial
bornology, the fine bornology. A subset $X \subset V$ belongs
to this bornology if and only $X$ is a bounded (in the usual
sense) subset of some finite dimensional subspace of $V$.
In this case $V$ is bornologically complete and any linear
map from $V$ to another bornological
vector space is bounded. By default we equip vector spaces
with a countable basis with the fine bornology.

More interestingly, if $V$ is a complete topological vector
space (e.g. a Fr\'echet space) we can define the precompact
bornology on $V$ as follows. We call $X \subset V$ bounded
if and only if its closure $\overline X$ is compact. Under
these assumptions $V$ is bornologically complete and
any continuous map between such vector spaces is bounded.
Conversely, any bounded linear map between two Fr\'echet
spaces with the precompact bornology is continuous
\cite[Lemma 2.2]{Mey0}.

The category of bornological vector spaces is not abelian,
but it does have enough injective and projective objects.
It also has inductive and projective limits.

Let $V$ be a bornological vector space and $\mr{End}_{bor} (V)$
the algebra of bounded linear maps $V \to V$. A subset
$L \subset \mr{End}_{bor} (V)$ is equibounded if
$L (X) := \{ l (x) : l \in L , x \in X \}$ is bounded for any
bounded set $X \subset V$. This gives $\mr{End}_{bor} (V)$ the
structure of a bornological algebra.

Let $A$ be a unital bornological algebra. By definition a
bornological $A$-module structure on $V$ is the same as a bounded
bilinear map $A \times V \to V$ or a bounded algebra homomorphism
$A \to \mr{End}_{bor} (V)$. Let $\mr{Mod}_{bor} (A)$ be the
category of bornological $A$-modules.

The $A$-balanced completed
bornological tensor product $\hot_A$ is defined by the following
universal property. Bounded linear maps $V_1 \hot_A V_2 \to V_3$
with $V_3$ complete correspond bijectively to bounded bilinear
maps $b : V_1 \times V_2 \to V_3$ that satisfy
$b (v_1 a, v_2 ) = b (v_1 , a v_2 )$.

In case $V_1 , V_2$ and $A$ have the fine bornology this is just
the algebraic tensor product over $A$. On the other hand, if
$V_1 ,V_2$ and $A$ are Fr\'echet spaces with the precompact
bornology then this agrees with the completed projective tensor
product over $A$.

By definition a sequence
\[
0 \to V_1 \to V_2 \to V_3 \to 0
\]
in $\mr{Mod}_{bor} (A)$ is a bornological extension if the maps are
bounded $A$-module homomorphisms and the sequence is split exact
in the category of bornological vector spaces.
We call a differential complex of bornological $A$-modules exact
if it admits a bounded $\mh C$-linear contraction. These notions
of extensions and exactness make $\mr{Mod}_{bor} (A)$ into an exact
category, whose derived category we denote by $\mr{Der}_{bor} (A)$.
Let $\hot_A^{\mh L}$ and $\mh R \mr{Hom}_A$ denote the total
derived functors of $\hot_A$ and $\mr{Hom}_A$. Thus
$U \hot_A^{\mh L} V$ is an object of $\mr{Der}_{bor} (A)$ whose
homology is $\mr{Tor}_*^A (U,V)$, and the (co)homology of
$\mh R \mr{Hom}_A (U,V)$ is $\mr{Ext}_A^* (U,V)$. However, the
total derived functors contain somewhat more information, as
the passage to homology forgets the bornological properties of
these differential complexes.

Suppose that $A, U$ and $V$ have the fine bornology. Then the
bornological functors $\hot_A$ and $\mr{Hom}_A$ agree with their
algebraic counterparts. Hence $\mr{Tor}^A_n (U,V)$ and $\mr{Ext}_A^n
(U,V)$ are the same in the algebraic and the bornological sense.

Let $f : A \to B$ be a morphism of unital complete
bornological algebras and
\[
0 \leftarrow A \leftarrow P_0 \leftarrow P_1 \leftarrow \cdots
\]
a resolution of $A$ by projective $A \hot A^{op}$-modules.

\begin{thm}\label{thm:2.15}
\cite[Theorem 35]{Mey1} The following are equivalent:
\begin{description}
\item[1)] $B \hot_A P_* \hot_A B$ is a projective
$B \hot B^{op}$-module resolution of $B$.
\item[2)] $(f^* B) \hot_A^{\mh L} (f^* B) \to B \hot_B^{\mh L} B
\; ( \cong B )$ is an isomorphism.
\item[3)] $(f^* U) \hot_A^{\mh L} (f^* V) \to U \hot_B^{\mh L} V$
is an isomorphism $\forall \, U \in \mr{Mod}_{bor} (B^{op}) , V \in
\mr{Mod}_{bor} (B)$.
\item[4)] $\mh R \mr{Hom}_B (U,V) \to \mh R \mr{Hom}_A (f^* U ,f^* V)$
is an isomorphism $\forall \, U,V \in \mr{Mod}_{bor} (B)$.
\item[5)] The functor $f^* : \mr{Der}_{bor} (B) \to \mr{Der}_{bor} (A)$
is fully faithful.
\end{description}
We call $f$ isocohomological if these conditions hold.
\end{thm}

Direct consequences of conditions 3) and 4) are
\begin{equation}\label{eq:2.31}
\begin{array}{lll}
\mr{Tor}_*^B (U,V) & \cong & \mr{Tor}_*^A (f^* U, f^* V) \vspace{1mm} \\
\mr{Ext}^*_B (U,V) & \cong & \mr{Ext}^*_A (f^* U, f^* V)
\end{array}
\end{equation}
where we mean are the derived functors in
the bornological category.
\\[1mm]

We equip $\mc H$ with the fine bornology and let $\mr{Mod}_{bor} (\mc H)$
be the category of all bornological $\mc H$-modules. Notice that any
$\mc H$-module can be made bornological by endowing it with the fine
bornology. This identifies Mod$ (\mc H)$ with a full subcategory of
$\mr{Mod}_{bor} (\mc H)$. An $\mc H$-module is bornologically projective
if and only if it is algebraically projective, namely if and only if
it is a direct summand of an (algebraically) free $\mc H$-module.
So as long as we are working in a purely algebraic setting the
bornological structure does not give much extra, but neither is it a
restriction.

We endow $\mc S$ with the precompact bornology, so that any finite
dimensional $\mc S$-module is bornological. We denote the category of
all bornological $\mc S$-modules by $\mr{Mod}_{bor} (\mc S)$. Probably
there exist $\mc S$-modules that do not admit the structure
of a bornological $\mc S$-module, but they seem to be
rather far-fetched. We note that a projective object of
$\mr{Mod}_{bor} (\mc S)$ is usually not a projective $\mc S$-module
in the algebraic sense, rather a completion of the latter.

A bornological $\mc H$-module $(\pi ,V)$ is called tempered if it
extends to $\mc S$, i.e. if the following equivalent conditions hold:
\begin{description}
\item[1)] $\pi$ extends to a bounded algebra homomorphism
$\mc S \to \mr{End}_{bor} (V)$
\item[2)] $\pi$ induces a bounded bilinear map $\mc S \times V \to V$
\end{description}

A (sub-)linear functional $f : \mc H \to \mh C$ is tempered if there
exist $C,N \in (0,\infty)$ such that
\[
|f (N_w)| \leq C (1 + \mc N (w) )^N \qquad \forall w \in W
\]
The collection of all tempered linear functionals is the continuous
dual space of $\mc S (\mc R ,q)$.

\begin{prop}\label{prop:2.14}
Let $V$ be a Fr\'echet space endowed with the precompact bornology.
An $\mc H$-module $(\pi ,V)$ is bornological if and only if
$\pi (h) : V \to V$ is continuous $\forall h \in \mc H$. Moreover it
is tempered if and only if the following equivalent conditions hold.
\begin{description}
\item[3)] $\pi$ induces a jointly continuous map
$\mc S \times V \to V$
\item[4)] $\pi$ induces a separately continuous map
$\mc S \times V \to V$
\item[5)] for every $v \in V$ and every continuous seminorm $p$ on $V$
the sublinear functional
\[
\mc H \to [0,\infty ) : h \to p (\pi (h) v)
\]
is tempered
\item[6)] for every $v \in V$ and every $f \in V^*$ the linear functional
\[
\mc H \to \mh C : h \to f (\pi (h) v)
\]
is tempered
\end{description}
In particular the category $\mr{Mod}_{Fr\acute e} (\mc S )$ of continuous
Fr\'echet $\mc S$-modules is a full subcategory of
$\mr{Mod}_{bor} (\mc S )$.
\end{prop}
\emph{Proof.} We already noted that $\pi (h) : V \to V$ is continuous
if and only if it is bounded. Since $\mc H$ carries the fine bornology
this is equivalent to the first assertion.

For the same reason $\mr{Mod}_{Fr\acute e} (\mc S )$ forms a full
subcategory of $\mr{Mod}_{bor} (\mc S )$.

It is clear that condition 3) implies the other five.
Conversely 3) follows from 2) by \cite[Lemma 2.2]{Mey0} and from 4)
by the Banach-Steinhaus theorem.

If $f \in V^*$ then $|f|$ is a continuous seminorm on $V$,
so 5) implies 6).

Finally we show that 6) implies 4). Endow $\mc H$ with the induced
topology from $\mc S$ and fix $v \in V$. By assumption the linear map
\begin{equation}\label{eq:2.30}
\mc H \to V : h \to \pi (h) v
\end{equation}
is continuous for the weak topology on $V$. Since $V$ is Fr\'echet
\eqref{eq:2.30} is also continuous for the metric topology on $V$
\cite[21.4.i]{KeNa}. Hence \eqref{eq:2.30} extends continuously to
the metric completion $\mc S$ of $\mc H$.

Now we fix $h = \sum\limits_{w \in W} h_w N_w \in \mc S$ and we write
$h_n = \sum\limits_{w : \mc N (w) \leq n} h_w N_w$. We assumed that $V$
is a Fr\'echet $\mc H$-module, so $(\pi (h_n ) )_{n=1}^\infty$ is
a sequence of continuous linear operators on $V$. We just showed that
for fixed $v \in V$ the sequence $(\pi (h_n) v )_{n=1}^\infty$
converges to $\pi (h) v$. The Banach-Steinhaus theorem (see e.g.
\cite[p. 104-105]{KeNa}) assures that $\pi (h)$ is continuous.

We conclude that $(h,v) \to \pi (h) v$ is separately continuous.
$\qquad \Box$ \\[3mm]

From the work of Casselman \cite[\S 4.4]{Cas} one can deduce more
concrete criteria for representations to be tempered or discrete series,
see \cite[Section 2.7]{Opd1}. It follows from these criteria that an
$\mc H$-module can only be tempered if all its $Z(W)$-weights
are unitary.